\documentclass[12pt]{article}
\usepackage[english]{babel}
\usepackage[utf8]{inputenc}
\usepackage[T1]{fontenc}
\usepackage{amsmath}
\usepackage{amsfonts}
\usepackage{amsthm}
\usepackage{amssymb}
\usepackage{graphicx}
\usepackage{dsfont}
\usepackage{color} 

\newtheorem{theorem}{Theorem}

\newtheorem{definition}[theorem]{Definition}
\newtheorem{lemma}[theorem]{Lemma}

\title{Numerical Solution for a Class of Evolution Differential Equations with $p$-Laplacian and Memory}\author{\textit{Jos\'{e} C.M. Duque\footnote{jduque@ubi.pt}, Belchior C.X. Mário\footnote{belchior.mario@ubi.pt} and Rui M.P. Almeida\footnote{ralmeida@ubi.pt}}\\ University of Beira Interior\\ Center of Mathematics and Applications\\ Covilh\~a, Portugal\\}
\date{\today}
\begin{document}

\maketitle

\begin{abstract}
In this paper we make a study of a partial integral differential equation with $p$-Laplacian using a mixed finite element method. Two stable and convergent fixed point schemes are proposed to solve the nonlinear algebraic system. Using the implementation of the method in Matlab environment, we numerically analyse the convergence with an example. Some other examples are presented in order to illustrate several asymptotic behaviours and some localization effects of the solutions.
\end{abstract}


\section{Introduction}

Partial differential equations are very comum to model various systems in science and engineer. But in some cases of heat transfer \cite{MR515894}, gas diffusion \cite{MR551044}, nuclear reactor dynamics \cite{MR706380} or mathematical finance \cite{MR2629564} it is needed to reflect the effects of the past history of the system in the model. A way to incorporate this effect in the model was found to be the inclusion of an integral term in the basic partial differential equation that leads to a Partial Integral Differential Equation (PIDE).

The $p$-Laplacian operator also appears in various applications of interest in many branches of science and engineering, such as flows in porous media \cite{MR2447796}, power-law materials \cite{MR760209}, nonlinear diffusion and filtration \cite{MR140343}, non-Newtonian fluids \cite{MR2221163} or elastic-plastic torsional creep \cite{MR533510} as example. The intersection of this two fields motivates the study of PIDEs with $p$-Laplacian.

In the present work we study the evolutionary intogrodifferential equation with $p$-Laplacian and memory
\begin{eqnarray} \label{eq:30}
\begin{cases}
u_{t}(x,t)-\Delta_{p}u(x,t)=\int_{0}^tg(t-s)\Delta_{p}u(x,s)ds+f(x,t),\,(x,t)\in\Omega\times]0,T],\\
u(x,t)=0,\,(x,t)\in\partial\Omega\times]0,T],\\
u(x,0)=u_{0}(x),\,x\in\Omega,
\end{cases}
\end{eqnarray}
where $f$ and $g$ are given functions. PIDEs of the above type arise naturally in many applications, such as, for instance, nonlocal reactive flows in porous media, heat conduction in materials with memory, and non-Fickian flow of fluid in porous media. Due to the presence of nonlinearity, it is difficult to find analytical solutions to these problems. Therefore, in practice, numerical simulation becomes a crucial tool to study the dynamics of this type of equations. The importance of our work lies in the study of the interaction between the $p$-Laplacian and the non-linear memory term on the solutions.

Considerable attention has been devoted to solving PIDEs numerically, as example, we can refer some recent works. To find the approximate solution of a linear PIDE, Avazzadeh et al. \cite{MR2880735} applied radial basis functions and the finite difference method. In \cite{MR3729702} a collocation method based on linear Legendre multiwavelets is developed for numerical solutions of PIDEs. In \cite{MR3936253} is proposed a finite element method in space and a Crank-Nicolson method in time to analyse a linear PIDE. Using reconstruction operators they obtained optimal order for a posteriori error estimates. In \cite{MR3958333} are presented some discretization techniques based on finite elements and fully implicit time discretization for solving parabolic PIDEs with nonlinear memory. The convergence order was proved to be $\mathcal{O}(h^r+\Delta t)$ for polynomials of degree $r-1$.

Since the pioneer works of \cite{GM75} and  \cite{DiB_book} several works on numerical simulations for the $p$-Laplacian were carried out. Recently, in \cite{MR2383205} is studied an adaptive FEM for a PDE with $p$-Laplacian using piecewise linear continuous functions and showed a linear convergence of the error. The authors of \cite{MR3463051} propose a hybridisable discontinuous Galerkin method for the $p$-Laplacian using polynomials of degree $r>0$. The numerical experiments showed optimal order of convergence. In \cite{MR4197298} a new phase field model involving the $p$-Laplacian operator is proposed. Optimal error estimates and convergence rates for the finite element approximation are proved.

Antontsev and his coauthors studied problem (\ref{eq:30}) with a nonlinear source term $\Theta(x,t,u)$, substituting the equation by a system composed of a diffusion-reaction equation and an integral equation. The authors proved that for $\max\{1,\frac{2n}{n+2}\}<p<\infty,g,g\prime\in L^2(0,T)$ and $u_{0}\in W_{0}^{1,p}(\Omega),f\in L^2(Q)$, the problem admits a weak solution that is local or global in time depending on the growth rate of $\Theta(x,t,s)$, when $\left|s\right|\rightarrow\infty$. They stated uniqueness conditions and they also proved that for $p>2$ and $s\Theta(x,t,s)\leq 0$, the data disturbances propagate with finite speed and the waiting time effect is possible \cite{MR3462616}.

In \cite{almeida2022mixed}, Almeida at al., presented a new mixed finite element method for integrodifferential equation with $p$-Lapacian and memory. In that work the existence, uniqueness and regularity of the discrete solutions are established. Error bounds, depending on the parameter $p$, are also obtained, however no simulations were done.

This paper is a continuation of the work carried out by Almeida at al. \cite{almeida2022mixed}, where some algorithms to solve the algebraic equations and several simulations in a matlab environment are presented, which verify and illustrate the theory developed in \cite{almeida2022mixed} and \cite{MR3462616}.

The present work has the following structure. This paper started with an introduction where, some themes related to the $p$-Laplacian theory, the memory term and the simulations are addressed. In section 2, we present the problem, some auxiliar theorems and lemmas and the main results already obtained. In section 3, two fixed point method schemes are proposed to solve the nonlinear system equations and it is proved their stability and convergence. In section 4, simulations are carried out, first to analyze the methods convergence and later to study its properties. Finally, section 5 ends the paper with the conclusions of the work carried out.

\section{Parabolic equation with $p$-Laplacian}

Let us consider the evolutionary integrodifferential equation with the homogeneous Dirichlet condition
\begin{eqnarray} \label{eq:14}
\begin{cases}
u_{t}-\Delta_{p}u=\int_{0}^{t}g(t-s)\Delta_{p}u(x,s)ds+f(x,t),\,(x,t)\in Q=\Omega\;\times]0,T],\\
u(x,t)=0,\,(x,t)\in\partial\Omega\times[0,T],\\
u(x,0)=u_{0}(x),\,x\in\Omega,
\end{cases}
\end{eqnarray}
where $u_0$, $g$ and $f$ are given functions, $\Omega\subset\mathbb{R}^n$ is a bounded domain with Lipschitz-continuous boundary. The $p$-Laplacian $\Delta_{p}u$ is given by
\begin{equation*}
\Delta_{p}u=\mathrm{div}\left(\left|\nabla u\right|^{p-2}\nabla u\right),\quad 1<p<\infty
\end{equation*}
and let us denote by
\begin{equation*} \label{eq:oc}
y(x,t)=\int_{0}^{t}g(t-s)\Delta_{p}u(x,s)ds
\end{equation*}
the memory term of the evolutionary integrodifferential equation.

Assuming
\begin{equation*}
g,g\prime\in L^2(0,T),\quad f\in L^2(Q),\quad u_0\in L^2(\Omega)\cap W_0^{1,p}(\Omega),
\end{equation*}
it is proved in \cite{MR3462616} that problema (\ref{eq:14}) has a unique weak solution.

In this work we will use the standard notations for norms and spaces, for more details and notation we refer to \cite{MR1801735,MR1625845,MR0259693}. Below we present an important lemma that can be found in \cite{MR1276708,MR3250760,Cho89}.
\begin{lemma} \label{eq:le1}
For every $p>1$ and $\delta\geq 0$ there are two positive constants $C_1$ and $C_2$, such that $\forall\zeta,\gamma\in\mathbb{R}^n$, $\zeta\neq\gamma$:
\begin{enumerate}
\item
\begin{equation*}
\left|\left|\zeta\right|^{p-2}\zeta-\left|\gamma\right|^{p-2}\gamma\right|\leq C_1\left|\zeta-\gamma \right|^{1-\delta}\left(\left|\zeta\right|+\left|\gamma\right|\right)^{p-2+\delta};
\end{equation*}
\item if $p>2$,
\begin{equation*}
\left(\left|\zeta\right|^{p-2}\zeta-\left|\gamma\right|^{p-2}\gamma,\zeta-\gamma\right)_{\mathbb{R}^n}\geq C_2\left|\zeta-\gamma\right|^{p}.
\end{equation*}
\end{enumerate}
\end{lemma}

\subsection{Auxiliary problem}

In the work developed by Almeida et al. in \cite{almeida2022mixed}, it was proved that the memory term satisfies the integral equation
\begin{eqnarray*} \label{eq:jaa}
y(x,t)=-\int_{0}^{t}g(t-s)y(x,s)ds+f_2(x,t,u),
\end{eqnarray*}
where $f_2$ is the nonlinear nonlocal operator
\begin{eqnarray*} \label{eq:f24}
f_2(x,t,u(x,t))&=&u(x,t)g(0)-u_0(x)g(t)+\int_{0}^{t}g\prime(t-s)u(x,s)ds\nonumber\\&&-\int_{0}^{t}g(t-s)f(x,s)ds.
\end{eqnarray*}

This allow us to consider the equivalent auxiliary problem of finding the pair $(u,y)$ that satisfies the conditions
\begin{eqnarray} \label{eq:ff23}
\begin{cases}
u_{t}-\Delta_{p}u=y(x,t)+f(x,t),\,(x,t)\in Q,\\
y(x,t)=-\int_{0}^{t}g(t-s)y(x,s)ds+f_2(x,t,u),\,(x,t)\in Q,\\
u(x,t)=0,\,(x,t)\in\partial\Omega\times[0,T],\\
u(x,0)=u_{0}(x),\,x\in\Omega,\\
y(x,0)=0,\,x\in\Omega.
\end{cases}
\end{eqnarray}

\begin{definition}
We say that $(u,y)\in(H_0^1(\Omega))^2$ is a weak solution of the evolutionary integrodifferential problem (\ref{eq:ff23}), if:
\begin{enumerate}
\item $u\in C([0,T],L^2(\Omega))\cap H_0^1(\Omega)$, $u_t\in L^2(\Omega)$, $\nabla u\in L^p\big(Q)$, $y\in L^2(Q)$, $u(x,0)=u_0(x)$, $y(x,0)=0$, $x\in\Omega$;
\item for each test function, $(w,v)\in(H_0^1(\Omega))^2$, the equalities
\end{enumerate}
\begin{eqnarray} \label{eq:ff27}
\begin{cases}
\int_{\Omega}^{}u_twdx+\int_{\Omega}^{}\left|\nabla u\right|^{p-2}\nabla u\nabla wdx=\int_{\Omega}^{}ywdx+\int_{\Omega}^{}fwdx,\,\forall w\in H_0^1(\Omega),\\
\int_{\Omega}^{}yvdx=-\int_{\Omega}^{}v\int_{0}^{t}g(t-s)y(x,s)dsdx+\int_{\Omega}^{}f_2vdx,\,\forall v\in H_0^1(\Omega),
\end{cases}
\end{eqnarray}
are valid.
\end{definition}

\section{Discretization in space}

\subsection{Lagrangian bases}

Let us consider $\mathcal{T}_h=\{T_0,\cdots,T_m\}$ a regular partition of $\Omega$ in simplexes with parameter $h$, and the space $\mathcal{S}^h\subset H_0^1(\Omega)$ defined by
\begin{equation*}
\mathcal{S}^h=\{w\in C^0(\Omega):w(x)=0,\;x\in\partial\Omega,\;w(x)|_{T_k}\in\mathcal{P}_r(T_k),\;k=0,\cdots,m\},
\end{equation*}
where $\mathcal{P}_r(T_k)$ is the set of polynomials of degree less than or equal to $r$ defined in $T_k$.

Below we find an important lemma where it is proved some inverse inequalities \cite{MR2050138}.
\begin{lemma}[Inverse estimates] \label{eq:le2}
Let $u^h\in\mathcal{S}^h$ and $n=1$, then there is a constant $C$ independent of $h$ such that
\begin{equation*}
\|\nabla u^h\|_{L^2(T_k)}\leq Ch^{-1}\|u^h\|_{L^2(T_k)}
\end{equation*}
and
\begin{equation*}
\|\nabla u^h\|_{L^\infty(T_k)}\leq Ch^{-\frac{3}{2}}\|u^h\|_{L^2(T_k)},\quad T_k\in\mathcal{T}_h.
\end{equation*}
\end{lemma}

The semi-discrete problem is to find $(u^h,y^h)\in \big(\mathcal{S}^h\big)^2$, such that
\begin{eqnarray} \label{eq:ff28}
\begin{cases}
\int_{\Omega}^{}u_t^hw^hdx+\int_{\Omega}^{}\left|\nabla u^h\right|^{p-2}\nabla u^h\nabla w^hdx=\int_{\Omega}^{}y^hw^hdx+\int_{\Omega}^{}fw^hdx,\,\forall w^h\in\mathcal{S}^h,\\
\int_{\Omega}^{}y^hv^hdx=-\int_{\Omega}^{}v^h\int_{0}^{t}g(t-s)y^h(x,s)dsdx+\int_{\Omega}^{}f_2^hv^hdx,\,\forall v^h\in\mathcal{S}^h,
\end{cases}
\end{eqnarray}
and
\begin{eqnarray*}
&&u^h(x,0)=u_0^h=\Pi_hu_0,\;\, y^h(x,0)=0,\quad\forall x\in\Omega,\\&& u^h(x,t)=0,\;\, y^h(x,t)=0,\quad\forall(x,t)\in\partial\Omega\times[0,T],
\end{eqnarray*}
where we will denote by $\Pi_h$ the interpolation operator into $\mathcal{S}^h$.

The following two theorems are proved in \cite{almeida2022mixed}.
\begin{theorem}
If $g,g\prime\in L^\infty(0,T)$, then there is a unique  solution of the semi-discrete problem (\ref{eq:ff28}) that satisfies
\begin{equation*}
\|u^h\|_{L^2(\Omega)}^2+\|\nabla u^h\|_{L^p(Q)}^p\leq C\|u_0\|_{L^2(\Omega)}^2+C\|f\|_{L^2(Q)}^2
\end{equation*}
and
\begin{equation*}
\|y^h\|_{L^2(\Omega)}^2\leq C\|u_0\|_{L^2(\Omega)}^2+C\|f\|_{L^2(Q)}^2,
\end{equation*}
where
\begin{equation*}
C=C\big(T,\|g\|_{L^\infty(0,T)},\|g\prime\|_{L^\infty(0,T)}\big).
\end{equation*}
\end{theorem}

\begin{theorem}[Convergence]
Let $(u,y)$ and $(u^h,y^h)$ be solutions of problems (\ref{eq:ff27}) and (\ref{eq:ff28}), respectively. If $g,g\prime\in L^\infty(0,T)$, $u_0\in H^{r+1}(\Omega)$, $f\in L^2(Q)$, then $\forall t\in [0,T]$,
\begin{eqnarray*} \label{eq:ff40}
\|u-u^h\|_{L^2(\Omega)}\leq Ch^{\frac{rp}{2(p-1)}}
\end{eqnarray*}
and
\begin{eqnarray*} \label{eq:ff41}
\|y-y^h\|_{L^2(\Omega)}\leq Ch^{\frac{rp}{2(p-1)}},
\end{eqnarray*}
where the constant $C$ does not depend on $h$ but may depend on $g$, $u$, $y$ and their derivatives.
\end{theorem}

\section{Discretization in time}

Following Almeida at al. \cite{almeida2022mixed}, we consider the partition $0=t_0<t_1<\cdots <t_N=T$, with step $\delta=\frac{T}{N}$, of $[0,T]$ and the notations
\begin{eqnarray*} \label{eq:ja1}
\bar{\partial}u^{\big(k+\frac{1}{2}\big)}&=&\frac{u^h(x,t_{k+1})-u^h(x,t_{k})}{\delta},\quad
\bar{u}^{\big(k+\frac{1}{2}\big)}=\frac{u^h(x,t_{k+1})+u^h(x,t_{k})}{2},\nonumber\\\bar{y}^{\big(k+\frac{1}{2}\big)}&=&\frac{y^h(x,t_{k+1})+y^h(x,t_{k})}{2},\quad 0\leq k\leq N.
\end{eqnarray*}

Applying the Crank-Nicolson method and the Trapezoidal quadrature, the totally discrete problem is to find $\big(u^h(x,t_{k+1}),y^h(x,t_{k+1})\big)\in(\mathcal{S}^h)^2$, solution of
\begin{eqnarray} \label{eq:ff31}
&&\int_{\Omega}^{}\bar{\partial}u^{\big(k+\frac{1}{2}\big)}w^hdx+\int_{\Omega}^{}\left|\nabla\bar{u}^{\big(k+\frac{1}{2}\big)}\right|^{p-2}\nabla\bar{u}^{\big(k+\frac{1}{2}\big)}\nabla w^hdx=\int_{\Omega}^{}\bar{y}^{\big(k+\frac{1}{2}\big)}w^hdx\nonumber\\&&+\int_{\Omega}^{}f\big(x,t_{k+\frac{1}{2}}\big)w^hdx
\end{eqnarray}
and
\begin{eqnarray} \label{eq:ff42}
&&\int_{\Omega}^{}\bar{y}^{\big(k+\frac{1}{2}\big)}v^hdx\nonumber\\&&=g(0)\int_{\Omega}^{}\bar{u}^{\big(k+\frac{1}{2}\big)}v^hdx-g\big(t_{k+\frac{1}{2}}\big)\int_{\Omega}^{}u_0^hv^hdx-Q_g(y^h)+Q_{g\prime}(u^h)\nonumber\\&&-I(f),
\end{eqnarray}
where
\begin{equation} \label{eq:ff311}
I(f)=\int_{0}^{t_{k+\frac{1}{2}}}g\big(t_{k+\frac{1}{2}}-s\big)\int_{\Omega}^{}f(x,s)v^h(x)dxds,
\end{equation}
\begin{eqnarray*} \label{eq:ja2}
&&Q_{g\prime}(u^h)\nonumber\\&&=\frac{\delta}{2}g\prime\big(t_{k+\frac{1}{2}}\big)\int_{\Omega}^{}u_0^h(x)v^h(x)dx+\delta\sum\limits_{j=1}^{k-1}g\prime\big(t_{k+\frac{1}{2}}-t_j\big)\int_{\Omega}^{}u^h(x,t_j)v^h(x)dx\nonumber\\&&+\frac{3\delta}{4}g\prime\big(t_{k+\frac{1}{2}}-t_k\big)\int_{\Omega}^{}u^h(x,t_k)v^h(x)dx+\frac{\delta}{8}g\prime(0)\int_{\Omega}^{}u^h(x,t_k)v^h(x)dx\nonumber\\&&+\frac{\delta}{8}g\prime(0)\int_{\Omega}^{}u^h(x,t_{k+1})v^h(x)dx
\end{eqnarray*}
and
\begin{eqnarray*}\label{eq:ja3}
&&Q_{g}(y^h)\nonumber\\&&=\frac{\delta}{2}g\big(t_{k+\frac{1}{2}}\big)\int_{\Omega}^{}y_0^h(x)v^h(x)dx+\delta\sum\limits_{j=1}^{k-1}g\big(t_{k+\frac{1}{2}}-t_j\big)\int_{\Omega}^{}y^h(x,t_j)v^h(x)dx\nonumber\\&&+\frac{3\delta}{4}g\big(t_{k+\frac{1}{2}}-t_k\big)\int_{\Omega}^{}y^h(x,t_k)v^h(x)dx+\frac{\delta}{8}g(0)\int_{\Omega}^{}y^h(x,t_k)v^h(x)dx\nonumber\\&&+\frac{\delta}{8}g(0)\int_{\Omega}^{}y^h(x,t_{k+1})v^h(x)dx.
\end{eqnarray*}

The following two theorems are proved in \cite{almeida2022mixed}.
\begin{theorem} \label{the3}
If $g,g\prime\in L^\infty(0,T)$, $u_0\in L^2(\Omega)$, $f\in L^2(Q)$, then for every $k\geq 0$, the discrete problem (\ref{eq:ff31}), (\ref{eq:ff42}) has a unique solution that satisfies
\begin{eqnarray*}
\|u^h(x,t_{k+1})\|_{L^2(\Omega)}^2\leq C\|u_0\|_{L^2(\Omega)}^2+C\|f\|_{L^2(Q)}^2+C\delta\sum\limits_{j=0}^{k}\|f\big(x,t_{j+\frac{1}{2}}\big)\|_{L^2(\Omega)}^2
\end{eqnarray*}
and
\begin{eqnarray*}
\|y^h(x,t_{k+1})\|_{L^2(\Omega)}^2\leq C\|u_0\|_{L^2(\Omega)}^2+C\|f\|_{L^2(Q)}^2+C\delta\sum\limits_{j=0}^{k}\|f\big(x,t_{j+\frac{1}{2}}\big)\|_{L^2(\Omega)}^2,
\end{eqnarray*}
where
\begin{equation*}
C=C\big(T,\|g\|_{L^\infty(0,T)},\|g\prime\|_{L^\infty(0,T)}\big).
\end{equation*}
\end{theorem}

\begin{theorem}[Convergence] \label{the4}
Let $(u,y)$ and $(u^h,y^h)$ be solutions to problems (\ref{eq:ff27}) and (\ref{eq:ff31})-(\ref{eq:ff42}), respectively. Suppose $g,g\prime\in L^\infty(0,T)$ and $u_0\in H^{r+1}(\Omega)$. If $\delta$ is small then for $0\leq k\leq N$
\begin{equation*}
\|u(x,t_k)-u^h(x,t_{k})\|_{L^2(\Omega)}\leq C\Big(h^{r+1}+\delta^2+h^{\frac{rp}{2(p-1)}}+\delta^{\frac{p}{p-1}}\Big)
\end{equation*}
and
\begin{equation*}
\|y(x,t_k)-y^h(x,t_{k})\|_{L^2(\Omega)}\leq C\Big(h^{r+1}+\delta^2+h^{\frac{rp}{2(p-1)}}+\delta^{\frac{p}{p-1}}\Big),
\end{equation*}
where the constant $C$ does not depend on $h$ or $\delta$ but may depend on $g$, $u$, $y$ and their derivatives.
\end{theorem}

Since we will only deal with functions of $\mathcal{S}^h$ henceforward we will omit the super index $h$. In order to simplify the notation we will use the super index $(j)$ in a function to represent the function evaluated at $t=t_j$.

In the numerical simulations we will consider the unidimensional case because is better to graphically illustrate the behaviour of the solution.

Rewriting the equation (\ref{eq:ff31}) of the form
\begin{eqnarray} \label{eq:j1}
&&2\int_{\Omega}^{}u^{(k+1)}wdx+\delta\int_{\Omega}{}\left|\frac{\nabla u^{(k+1)}+\nabla u^{(k)}}{2}\right|^{p-2}\left(\nabla u^{(k+1)}-\nabla u^{(k)}\right)\nabla wdx\nonumber\\&&=\delta\int_{\Omega}{} \big(y^{(k+1)}-y^{(k)}\big)wdx+2\delta\int_{\Omega}{}f^{\big(k+\frac{1}{2}\big)}wdx\nonumber\\&&+2\int_{\Omega}{}u^{(k)}wdx,\;\forall w\in\mathcal{S}^h
\end{eqnarray}
and the equation (\ref{eq:ff42}) of the form
\begin{eqnarray} \label{eq:j2}
&&\left(1+\frac{\delta}{4}g(0)\right)\int_{\Omega}{}y^{(k+1)}vdx-\left(2g(0)+\frac{\delta}{4}g\prime(0)\right)\int_{\Omega}^{}u^{(k+1)}vdx\nonumber\\&&=\left(-\frac{1}{2}+\frac{3\delta}{4}g\big(t_{k+\frac{1}{2}}-t_k\big)+\frac{\delta}{8}g(0)\right)\int_{\Omega}^{}y^{(k)}vdx\nonumber\\&&+\delta\sum\limits_{l=1}^{k-1}g\left(t_{k+\frac{1}{2}}-t_l\right)\int_{\Omega}^{}y^{(l)}vdx+\frac{\delta}{2}g\big(t_{k+\frac{1}{2}}\big)\int_{\Omega}^{}y^{(0)}vdx\nonumber\\&&+\left(\frac{1}{2}g(0)-\frac{3\delta}{4}g\prime\big(t_{k+\frac{1}{2}}-t_k\big)-\frac{\delta}{8}g(0)\right)\int_{\Omega}^{}u^{(k)}vdx\nonumber\\&&-\left(g\big(t_{k+\frac{1}{2}}\big)+\frac{\delta}{2}g\prime\big(t_{k+\frac{1}{2}}\big)\right)\int_{\Omega}^{}u^{(0)}vdx+\delta\sum\limits_{l=1}^{k-1}g\prime\left(t_{k+\frac{1}{2}}-t_l\right)\int_{\Omega}^{}u^{(l)}vdx\nonumber\\&&+\frac{\delta}{4}g\big(t_{k+\frac{1}{2}}\big)\int_{\Omega}^{}f^{(0)}vdx+\frac{3\delta}{4}g\left(t_{k+\frac{1}{2}}-t_{\frac{1}{2}}\right)\int_{\Omega}^{}f^{\big(\frac{1}{2}\big)}vdx\nonumber\\&&+\delta\sum\limits_{m=1}^{k}g\left(t_{k+\frac{1}{2}}-t_{m+\frac{1}{2}}\right)\int_{\Omega}^{}f^{\big(m+\frac{1}{2}\big)}vdx\nonumber\\&&+\frac{\delta}{2}g(0)\int_{\Omega}^{}f^{\big(k+\frac{1}{2}\big)}vdx,\;\forall v\in \mathcal{S}^h,
\end{eqnarray}
we obtain a system of nonlinear algebraic equations with the unknown\linebreak $\big(u^{(k+1)},y^{(k+1)}\big)$. To solve it in the case $p\geq 3$, we propose the fixed point method with the following iterative scheme: Given $f,g,h,\delta,u^{(0)},\cdots,u^{(k)},\linebreak y^{(0)},\cdots,y^{(k)} $, we consider $u^{(k+1)}$ and $y^{(k+1)}$ the limits of the sequences $\big(\mathcal{U}_{(n) }\big)$ and $\big(\mathcal{Y}_{(n)}\big)$ defined by
\begin{eqnarray}\label{eq:j3}
&&2\int_{\Omega}{}\mathcal{U}_{(n+1)}wdx+\delta\int_{\Omega}{}\left|\frac{\nabla\mathcal{U}_{(n)}+\nabla u^{(k)}}{2}\right|^{p-2}\left(\nabla\mathcal{U}_{(n+1)}+\nabla u^{(k)}\right)\nabla wdx\nonumber\\&&=\delta\int_{\Omega}{}\big(\mathcal{Y}_{(n+1)}-y^{(k)}\big)wdx+2\delta\int_{\Omega}{}f^{\big(k+\frac{1}{2}\big)}wdx\nonumber\\&&+2\int_{\Omega}{}u^{(k)}wdx,\;\forall w\in \mathcal{S}^h
\end{eqnarray}
and
\begin{eqnarray}\label{eq:j4}
&&\left(1+\frac{\delta}{4}g(0)\right)\int_{\Omega}{}\mathcal{Y}_{(n+1)}vdx-\left(2g(0)+\frac{\delta}{4}g\prime(0)\right)\int_{\Omega}{}\mathcal{U}_{(n+1)}vdx\nonumber\\&&=\left(-\frac{1}{2}+\frac{3\delta}{4}g\big(t_{k+\frac{1}{2}}-t_k\big)+\frac{\delta}{8}g(0)\right)\int_{\Omega}^{}y^{(k)}vdx\nonumber\\&&+\delta\sum\limits_{l=1}^{k-1}g\left(t_{k+\frac{1}{2}}-t_l\right)\int_{\Omega}^{}y^{(l)}vdx+\frac{\delta}{2}g\big(t_{k+\frac{1}{2}}\big)\int_{\Omega}^{}y^{(0)}vdx\nonumber\\&&+\left(\frac{1}{2}g(0)-\frac{3\delta}{4}g\prime\big(t_{k+\frac{1}{2}}-t_k\big)-\frac{\delta}{8}g(0)\right)\int_{\Omega}^{}u^{(k)}vdx\nonumber\\&&-\left(g\big(t_{k+\frac{1}{2}}\big)+\frac{\delta}{2}g\prime\big(t_{k+\frac{1}{2}}\big)\right)\int_{\Omega}^{}u^{(0)}vdx+\delta\sum\limits_{l=1}^{k-1}g\prime\left(t_{k+\frac{1}{2}}-t_l\right)\int_{\Omega}^{}u^{(l)}vdx\nonumber\\&&+\frac{\delta}{4}g\big(t_{k+\frac{1}{2}}\big)\int_{\Omega}^{}f^{(0)}vdx+\frac{3\delta}{4}g\left(t_{k+\frac{1}{2}}-t_{\frac{1}{2}}\right)\int_{\Omega}^{}f^{\big(\frac{1}{2}\big)}vdx\nonumber\\&&+\delta\sum\limits_{m=1}^{k}g\left(t_{k+\frac{1}{2}}-t_{m+\frac{1}{2}}\right)\int_{\Omega}^{}f^{\big(m+\frac{1}{2}\big)}vdx\nonumber\\&&+\frac{\delta}{2}g(0)\int_{\Omega}^{}f^{\big(k+\frac{1}{2}\big)}vdx,\;\forall v\in \mathcal{S}^h,
\end{eqnarray}
where
\begin{equation} \label{eq:d1}
\mathcal{U}_{(0)}=u^{(k)}\quad\mathrm{and}\quad\mathcal{Y}_{(0)}=y^{(k)}.
\end{equation}

The system (\ref{eq:j3})-(\ref{eq:j4}) is linear that has a unique solution.

Next Lemma shows that the sequences are stable.
\begin{lemma} \label{eq:led}
Let $\mathcal{U}_{(n+1)}$, $\mathcal{Y}_{(n+1)}$ be solutions of (\ref{eq:j3}) and (\ref{eq:j4}). If $f\in L^2(Q),\linebreak g,g\prime\in L^\infty(0,T)$, $u_0\in L^2(\Omega)$ and $\delta$ is small. Then
\begin{equation} \label{eq:d2}
\|\mathcal{U}_{(n+1)}\|^2_{L^2(\Omega)}<C
\end{equation}
and
\begin{equation} \label{eq:d3}
\|\mathcal{Y}_{(n+1)}\|^2_{L^2(\Omega)}<C,\quad n\in\mathbb{N},
\end{equation}
where
\begin{equation*}
C=C\left(\|f\|^2_{L^2(Q)},\|g\|^2_{L^\infty(0,T)},\|g\prime\|^2_{L^2(0,T)},\|u_{0}\|^2_{L^2(\Omega)}\right).
\end{equation*}
\end{lemma}
\begin{proof}
For $n=0$, the conditions (\ref{eq:d2}) and (\ref{eq:d3}) are true by condition (\ref{eq:d1}) and Theorem \ref{the3}. Suppose now that $\mathcal{U}_{(l)}$ and $\mathcal{Y}_{(l)}$ satisfies (\ref{eq:d2}) and (\ref{eq:d3}), for $l=0,\cdots,n$. If we consider in (\ref{eq:j4}) $v=\mathcal{Y}_{(n+1)}$ then
\begin{eqnarray} \label{eq:j5}
&&\left(1+\frac{\delta}{4}g(0)\right)\int_{\Omega}{}\mathcal{Y}_{(n+1)}\mathcal{Y}_{(n+1)}dx
\nonumber\\&&=\left(2g(0)+\frac{\delta}{4}g\prime(0)\right)\int_{\Omega}{}\mathcal{U}_{(n+1)}\mathcal{Y}_{(n+1)}dx+\left(-\frac{1}{2}+\frac{3\delta}{4}g\big(t_{k+\frac{1}{2}}-t_k\big)\right.\nonumber\\&&+\left.\frac{\delta}{8}g(0)\right)\int_{\Omega}{}y^{(k)}\mathcal{Y}_{(n+1)}dx+\delta\sum\limits_{l=1}^{k-1}g\left(t_{k+\frac{1}{2}}-t_l\right)\int_{\Omega}{}y^{(l)}\mathcal{Y}_{(n+1)}dx\nonumber\\&&+\frac{\delta}{2}g\big(t_{k+\frac{1}{2}}\big)\int_{\Omega}{}y^{(0)}\mathcal{Y}_{(n+1)}dx+\left(\frac{1}{2}g(0)-\frac{3\delta}{4}g\prime\big(t_{k+\frac{1}{2}}-t_k\big)\right.\nonumber\\&&-\left.\frac{\delta}{8}g(0)\right)\int_{\Omega}{}u^{(k)}\mathcal{Y}_{(n+1)}dx-\left(g\big(t_{k+\frac{1}{2}}\big)+\frac{\delta}{2}g\prime\big(t_{k+\frac{1}{2}}\big)\right)\int_{\Omega}{}u^{(0)}\mathcal{Y}_{(n+1)}dx\nonumber\\&&+\delta\sum\limits_{l=1}^{k-1}g\prime\left(t_{k+\frac{1}{2}}-t_l\right)\int_{\Omega}{}u^{(l)}\mathcal{Y}_{(n+1)}dx+\frac{\delta}{4}g\big(t_{k+\frac{1}{2}}\big)\int_{\Omega}{}f^{(0)}\mathcal{Y}_{(n+1)}dx\nonumber\\&&+\frac{3\delta}{4}g\left(t_{k+\frac{1}{2}}-t_{\frac{1}{2}}\right)\int_{\Omega}{}f^{\big(\frac{1}{2}\big)}\mathcal{Y}_{(n+1)}dx+\delta\sum\limits_{m=1}^{k}g\left(t_{k+\frac{1}{2}}-t_{m+\frac{1}{2}}\right)\nonumber\\&&\times\int_{\Omega}{}f^{\big(m+\frac{1}{2}\big)}\mathcal{Y}_{(n+1)}dx+\frac{\delta}{2}g(0)\int_{\Omega}{}f^{\big(k+\frac{1}{2}\big)}\mathcal{Y}_{(n+1)}dx.
\end{eqnarray}
Applying Young's inequality to (\ref{eq:j5}) and Theorem \ref{the3}, we have
\begin{eqnarray*}
\left(1+\frac{\delta}{4}g(0)\right)\|\mathcal{Y}_{(n+1)}\|^2_{L^2(\Omega)}\leq C\|\mathcal{U}_{(n+1)}\|^2_{L^2(\Omega)}+\epsilon\|\mathcal{Y}_{(n+1)}\|^2_{L^2(\Omega)}+C.
\end{eqnarray*}
If $\delta$ is small enough (only for the case $g(0)<0$) and $\epsilon$ is adequate, we get
\begin{eqnarray}\label{eq:j6}
\|\mathcal{Y}_{(n+1)}\|^2_{L^2(\Omega)}\leq C\|\mathcal{U}_{(n+1)}\|^2_{L^2(\Omega)}+C.
\end{eqnarray}
Returning to the equation (\ref{eq:j3}), with $w=\mathcal{U}_{(n+1)}+u^{(k)}$, there is
\begin{eqnarray}\label{eq:j7}
&&2\|\mathcal{U}_{(n+1)}\|^2_{L^2(a,b)}+\delta\int_{\Omega}^{}\left|\frac{\nabla\mathcal{U}_{(n)}+\nabla u^{(k)}}{2}\right|^{p-2}\left(\nabla\mathcal{U}_{(n+1)}+\nabla u^{(k)}\right)^2dx\nonumber\\&&=\delta\int_{\Omega}{}\big(\mathcal{Y}_{(n+1)}+y^{(k)}\big)\big(\mathcal{U}_{(n+1)}+u^{(k)}\big)dx+2\delta\int_{\Omega}{}f^{\big(k+\frac{1}{2}\big)}\big(\mathcal{U}_{(n+1)}+u^{(k)}\big)dx\nonumber\\&&+\int_{\Omega}{}u^{(k)}\big(\mathcal{U}_{(n+1)}+u^{(k)}\big)dx-2\int_{\Omega}{}\mathcal{U}_{(n+1)}u^{(k)}dx.
\end{eqnarray}
Applying Young's inequality to (\ref{eq:j7}), we get
\begin{eqnarray*}
(1-C\delta^2)\|\mathcal{U}_{(n+1)}\|^2_{L^2(\Omega)}\leq C\delta^2+C.
\end{eqnarray*}
If $\delta$ is small enough such that $(1-C\delta^2)>0$, then
\begin{eqnarray*}
\|\mathcal{U}_{(n+1)}\|^2_{L^2(\Omega)}\leq C.
\end{eqnarray*}
Going back to $\mathcal{Y}_{(n+1)}$ with the equation (\ref{eq:j6}), we have
\begin{eqnarray*}
\|\mathcal{Y}_{(n+1)}\|^2_{L^2(\Omega)}\leq C,
\end{eqnarray*}
what proves the intended.
\end{proof}
\begin{theorem}[Convergence] \label{the1}
Under the conditions of Lemma \ref{eq:led}, if $\delta$ and $h$ are appropriate, then the successions $\big(\mathcal{U}_{(n)}\big)$ and $\big(\mathcal{Y}_{(n)}\big)$ defined by (\ref{eq:j3}) and (\ref{eq:j4}) are convergent.
\end{theorem}
\begin{proof}
Subtracting the equation (\ref{eq:j4}) in the iteration $(n)$ from the equation (\ref{eq:j4}) in the iteration $(n+1)$, we have
\begin{equation*}
\left(1+\frac{\delta}{4}g(0)\right)\int_{\Omega}{}\big(\mathcal{Y}_{(n+1)}-\mathcal{Y}_{(n)}\big)vdx-\left(2g(0)+\frac{\delta}{4}g\prime(0)\right)\int_{\Omega}{}\big(\mathcal{U}_{(n+1)}-\mathcal{U}_{(n)}\big)vdx=0.
\end{equation*}
Making $v=\mathcal{Y}_{(n+1)}-\mathcal{Y}_{(n)}$,
\begin{eqnarray}\label{eq:h1}
&&\left(1+\frac{\delta}{4}g(0)\right)\|\mathcal{Y}_{(n+1)}-\mathcal{Y}_{(n)}\|^2_{L^2(\Omega)}\nonumber\\&&=\left(2g(0)+\frac{\delta}{4}g\prime(0)\right)\int_{\Omega}{}\big(\mathcal{U}_{(n+1)}-\mathcal{U}_{(n)}\big)\big(\mathcal{Y}_{(n+1)}-\mathcal{Y}_{(n)}\big)dx.
\end{eqnarray}
Applying Young's inequality to (\ref{eq:h1}), we get
\begin{eqnarray*}\label{eq:h}
\|\mathcal{Y}_{(n+1)}-\mathcal{Y}_{(n)}\|^2_{L^2(\Omega)}\leq C\|\mathcal{U}_{(n+1)}-\mathcal{U}_{(n)}\|^2_{L^2(\Omega)}.
\end{eqnarray*}
Subtracting the equation (\ref{eq:j3}) in the iteration $(n)$ and the equation (\ref{eq:j3}) in the iteration $(n+1)$, there is
\begin{eqnarray*}
&&2\int_{\Omega}{}\big(\mathcal{U}_{(n+1)}-\mathcal{U}_{(n)}\big)wdx\nonumber\\&&+\delta\int_{\Omega}{}\left|\frac{\nabla\mathcal{U}_{(n)}+\nabla u^{(k)}}{2}\right|^{p-2}\left(\nabla\mathcal{U}_{(n+1)}-\nabla\mathcal{U}_{(n)}\right)\nabla wdx\nonumber\\&&=-\delta\int_{\Omega}{}\left(\left|\frac{\nabla\mathcal{U}_{(n)}+\nabla u^{(k)}}{2}\right|^{p-2}-\left|\frac{\nabla\mathcal{U}_{(n-1)}+\nabla u^{(k)}}{2}\right|^{p-2}\right)\left(\nabla\mathcal{U}_{(n)}\right.\nonumber\\&&+\left.\nabla u^{(k)}\right)\nabla wdx+\delta\int_{\Omega}{}\big(\mathcal{Y}_{(n+1)}-\mathcal{Y}_{(n)}\big)wdx.
\end{eqnarray*}
Considering $w=\mathcal{U}_{(n+1)}-\mathcal{U}_{(n)}$,
\begin{eqnarray}\label{eq:h2}
&&2\|\mathcal{U}_{(n+1)}-\mathcal{U}_{(n)}\|^2_{L^2(\Omega)}\nonumber\\&&+\delta\int_{\Omega}^{}\left|\frac{\nabla\mathcal{U}_{(n)}+\nabla u^{(k)}}{2}\right|^{p-2}\left(\nabla\mathcal{U}_{(n+1)}-\nabla\mathcal{U}_{(n)}\right)^2dx\nonumber\\&&=-\delta\int_{\Omega}{}\left(\left|\frac{\nabla\mathcal{U}_{(n)}+\nabla u^{(k)}}{2}\right|^{p-2}-\left|\frac{\nabla\mathcal{U}_{(n-1)}+\nabla u^{(k)}}{2}\right|^{p-2}\right)\nonumber\\&&\times\Big(\nabla\mathcal{U}_{(n)}+\nabla u^{(k)}\Big)\Big(\nabla\mathcal{U}_{(n+1)}-\nabla\mathcal{U}_{(n)}\Big)dx\nonumber\\&&+\delta\int_{\Omega}^{}\big(\mathcal{Y}_{(n+1)}-\mathcal{Y}_{(n)}\big)\big(\mathcal{U}_{(n+1)}-\mathcal{U}_{(n)}\big)dx.
\end{eqnarray}
Applying Young's inequality to (\ref{eq:h2}) and Lemma \ref{eq:le2}, we have
\begin{eqnarray*}
&&\big(1-C\delta^2-C\delta h^{-2}\big)\|\mathcal{U}_{(n+1)}-\mathcal{U}_{(n)}\|^2_{L^2(\Omega)}\nonumber\\&&\leq C\delta h^{-3}\left\|\left|\frac{\nabla\mathcal{U}_{(n)}+\nabla u^{(k)}}{2}\right|^{p-2}-\left|\frac{\nabla\mathcal{U}_{(n-1)}+\nabla u^{(k)}}{2}\right|^{p-2}\right\|^2_{L^2(\Omega)}.
\end{eqnarray*}
Applying Lemma \ref{eq:le2}, there is
\begin{eqnarray*}
\big(1-C\delta^2-C\delta h^{-2}\big)\|\mathcal{U}_{(n+1)}-\mathcal{U}_{(n)}\|^2_{L^2(\Omega)}\leq C\delta h^{6(1-p)}\|\mathcal{U}_{(n)}-\mathcal{U}_{(n-1)}\|^2_{L^2(\Omega)}.
\end{eqnarray*}
If $\delta$ and $h$ are small enough such that $\big(1-C\delta^2-C\delta h^{-2}\big)>0$ and\\ $\frac{C\delta h^{6(1-p)}}{1-C\delta^2-C\delta h^{-2}}<1$, we have
\begin{eqnarray}\label{eq:h3}
\|\mathcal{U}_{(n+1)}-\mathcal{U}_{(n)}\|^2_{L^2(\Omega)}\leq C\|\mathcal{U}_{(n)}-\mathcal{U}_{(n-1)}\|^2_{L^2(\Omega)},\quad \mathrm{with}\;\,0<C<1.
\end{eqnarray}
Iterating the equation (\ref{eq:h3}), we have
\begin{eqnarray*}
\|\mathcal{U}_{(n+1)}-\mathcal{U}_{(n)}\|^2_{L^2(\Omega)}\leq C^n\|\mathcal{U}_{(1)}-\mathcal{U}_{(0)}\|^2_{L^2(\Omega)}\rightarrow 0,n\rightarrow +\infty.
\end{eqnarray*}
Going back to the equation (\ref{eq:h}), we have
\begin{eqnarray*}
\|\mathcal{Y}_{(n+1)}-\mathcal{Y}_{(n)}\|^2_{L^2(\Omega)}\rightarrow 0,n\rightarrow +\infty,
\end{eqnarray*}
what proves the intended.
\end{proof}

To solve the system of nonlinear algebraic equations (\ref{eq:j1}) and (\ref{eq:j2}), for the case $2<p<3$, taking into account the use of the method of the fixed point, let's replace the equation (\ref{eq:j3}) with
\begin{eqnarray}\label{eq:h13}
&&2\int_{\Omega}^{}\mathcal{U}_{(n+1)}wdx+\delta\int_{\Omega}^{}\left|\frac{\nabla\mathcal{U}_{(n)}+\nabla u^{(k)}}{2}\right|^{p-2}\Big(\nabla\mathcal{U}_{(n)}+\nabla u^{(k)}\Big)\nabla wdx\nonumber\\&&=\delta\int_{\Omega}^{}\big(\mathcal{Y}_{(n+1)}-y^{(k)}\big)wdx+\int_{\Omega}^{}f^{\big(k+\frac{1}{2}\big)}wdx+\int_{\Omega}^{}u^{(k)}wdx.
\end{eqnarray}
The system (\ref{eq:h13})-(\ref{eq:j4}) is also a linear system and has a unique solution for each $n\in\mathbb{N}$.
\begin{lemma}\label{eq:led1}
Let $\mathcal{U}_{(n+1)}$, $\mathcal{Y}_{(n+1)}$ be solutions of (\ref{eq:j4}) and (\ref{eq:h13}). If $f\in L^2(Q),\linebreak g,g\prime\in L^\infty(0,T)$, $u_0\in L^2(\Omega)$ and $\delta$ is small. Then we have the conditions (\ref{eq:d2}) and (\ref{eq:d3}).
\end{lemma}
\begin{proof}
For $n=0$, the conditions (\ref{eq:d2}) and (\ref{eq:d3}) are true by condition (\ref{eq:d1}) and Theorem \ref{the3}. Suppose now that $\mathcal{U}_{(l)}$ and $\mathcal{Y}_{(l)}$ satisfies (\ref{eq:d2}) and (\ref{eq:d3}), for $l=0,\cdots,n$. If we consider in (\ref{eq:j4}) $v=\mathcal{Y}_{(n+1)}$ and argue like in Lemma \ref{eq:led} we get
\begin{eqnarray}\label{eq:h12}
\|\mathcal{Y}_{(n+1)}\|^2_{L^2(\Omega)}\leq C\|\mathcal{U}_{(n+1)}\|^2_{L^2(\Omega)}+C.
\end{eqnarray}
Considering equation (\ref{eq:h13}), with $w=\mathcal{U}_{(n+1)}$, we have
\begin{eqnarray}\label{eq:h6}
&&\|\mathcal{U}_{(n+1)}\|^2_{L^2(\Omega)}+\delta\int_{\Omega}^{}\left|\frac{\nabla\mathcal{U}_{(n)}+\nabla u^{(k)}}{2}\right|^{p-2}\Big(\nabla\mathcal{U}_{(n)}+\nabla u^{(k)}\Big)\nabla\mathcal{U}_{(n+1)}dx\nonumber\\&&=\delta\int_{\Omega}^{}\mathcal{Y}_{(n+1)}\mathcal{U}_{(n+1)}dx+\delta\int_{\Omega}^{}y^{(k)}\mathcal{U}_{(n+1)}dx+2\delta\int_{\Omega}^{}f^{\big(k+\frac{1}{2}\big)}\mathcal{U}_{(n+1)}dx.
\end{eqnarray}
Applying Young's inequality to (\ref{eq:h6}), we get
\begin{eqnarray*}
&&\left(\frac{1}{2}-C\delta-C\delta h^{-2}\right)\|\mathcal{U}_{(n+1)}\|^2_{L^2(\Omega)}\nonumber\\&&\leq C\delta^2\|y^{(k)}\|^2_{L^2(\Omega)}+C\delta^2\|f^{\big(k+\frac{1}{2}\big)}\|^2_{L^2(\Omega)}+C\delta h^{3(1-p)}\|\mathcal{U}_{(n)}\|^{p-1}_{L^{2}(\Omega)}+C.
\end{eqnarray*}
If $\delta$ and $h$ are suitable, then
\begin{eqnarray}\label{eq:h7}
\|\mathcal{U}_{(n+1)}\|^2_{L^2(\Omega)}\leq C\delta+C\delta h^{2(1-p)}\|\mathcal{U}_{(n)}\|^{p-1}_{L^{2}(\Omega)}.
\end{eqnarray}
By the hypothesis, we have
\begin{eqnarray*}
\|\mathcal{U}_{(n+1)}\|^2_{L^2(\Omega)}\leq C.
\end{eqnarray*}
Going back to $\mathcal{Y}_{(n+1)}$ with the equation (\ref{eq:h12}), we have
\begin{eqnarray*}
\|\mathcal{Y}_{(n+1)}\|^2_{L^2(\Omega)}\leq C,
\end{eqnarray*}
what proves the intended.
\end{proof}
\begin{theorem}[Convergence] \label{the2}
Under the conditions of Lemma \ref{eq:led1} if $\delta$ and $h$ are appropriate, then the successions $\big(\mathcal{U}_{(n)}\big)$ and $\big(\mathcal{Y}_{(n)}\big)$ defined by (\ref{eq:j4}) and (\ref{eq:h13}) are convergent.
\end{theorem}
\begin{proof}
Subtracting the equation (\ref{eq:j4}) in the iteration $(n)$ from the equation (\ref{eq:j4}) in the iteration $(n+1)$, we have
\begin{equation*}
\left(1+\frac{\delta}{4}g(0)\right)\int_{\Omega}^{}\big(\mathcal{Y}_{(n+1)}-\mathcal{Y}_{(n)}\big)vdx-\left(2g(0)+\frac{\delta}{4}g\prime(0)\right)\int_{\Omega}^{}\big(\mathcal{U}_{(n+1)}-\mathcal{U}_{(n)}\big)vdx=0.
\end{equation*}
Making $v=\mathcal{Y}_{(n+1)}-\mathcal{Y}_{(n)}$,
\begin{eqnarray}\label{eq:h14}
&&\left(1+\frac{\delta}{4}g(0)\right)\|\mathcal{Y}_{(n+1)}-\mathcal{Y}_{(n)}\|^2_{L^2(\Omega)}\nonumber\\&&=\left(2g(0)+\frac{\delta}{4}g\prime(0)\right)\int_{\Omega}^{}\big(\mathcal{U}_{(n+1)}-\mathcal{U}_{(n)}\big)\big(\mathcal{Y}_{(n+1)}-\mathcal{Y}_{(n)}\big)dx.
\end{eqnarray}
Applying Young's inequality to (\ref{eq:h14}), we get
\begin{eqnarray}\label{eq:h10}
\|\mathcal{Y}_{(n+1)}-\mathcal{Y}_{(n)}\|^2_{L^2(\Omega)}\leq C\|\mathcal{U}_{(n+1)}-\mathcal{U}_{(n)}\|^2_{L^2(\Omega)}.
\end{eqnarray}
Subtracting the equation (\ref{eq:h13}) in the iteration $(n)$ and the equation (\ref{eq:h13}) in the iteration $(n+1)$, there is
\begin{eqnarray*}
&&2\int_{\Omega}^{}\big(\mathcal{U}_{(n+1)}-\mathcal{U}_{(n)}\big)wdx+\delta\int_{\Omega}^{}\left(\left|\frac{\nabla\mathcal{U}_{(n)}+\nabla u^{(k)}}{2}\right|^{p-2}\left(\nabla\mathcal{U}_{(n)}+\nabla u^{(k)}\right)\right.\\&&\left.-\left|\frac{\nabla\mathcal{U}_{(n-1)}+\nabla u^{(k)}}{2}\right|^{p-2}\left(\nabla\mathcal{U}_{(n-1)}+\nabla u^{(k)}\right)\right)\nabla wdx\nonumber\\&&=\delta\int_{\Omega}^{}\big(\mathcal{Y}_{(n+1)}-\mathcal{Y}_{(n)}\big)wdx,
\end{eqnarray*}
considering $w=\mathcal{U}_{(n+1)}-\mathcal{U}_{(n)}$,
\begin{eqnarray}\label{eq:h15}
&&2\|\mathcal{U}_{(n+1)}-\mathcal{U}_{(n)}\|^2_{L^2(\Omega)}\nonumber\\&&=-\delta\int_{\Omega}^{}\left(\left|\frac{\nabla\mathcal{U}_{(n)}+\nabla u^{(k)}}{2}\right|^{p-2}\left(\nabla\mathcal{U}_{(n)}+\nabla u^{(k)}\right)\right.\nonumber\\&&\left.-\left|\frac{\nabla\mathcal{U}_{(n-1)}+\nabla u^{(k)}}{2}\right|^{p-2}\left(\nabla\mathcal{U}_{(n-1)}+\nabla u^{(k)}\right)\left(\nabla\mathcal{U}_{(n+1)}-\nabla\mathcal{U}_{(n)}\right)\right)dx\nonumber\\&&+\delta\int_{\Omega}{}\big(\mathcal{Y}_{(n+1)}-\mathcal{Y}_{(n)}\big)\big(\mathcal{U}_{(n+1)}-\mathcal{U}_{(n)}\big)dx.
\end{eqnarray}
Applying Young's inequality to (\ref{eq:h15}), we have
\begin{eqnarray}\label{eq:h16}
&&2\|\mathcal{U}_{(n+1)}-\mathcal{U}_{(n)}\|^2_{L^2(\Omega)}\nonumber\\&&\leq C\delta\left\|\left|\frac{\nabla\mathcal{U}_{(n)}+\nabla u^{(k)}}{2}\right|^{p-2}\left(\nabla\mathcal{U}_{(n)}+\nabla u^{(k)}\right)\right.\nonumber\\&&\left.-\left|\frac{\nabla\mathcal{U}_{(n-1)}+\nabla u^{(k)}}{2}\right|^{p-2}\left(\nabla\mathcal{U}_{(n-1)}+\nabla u^{(k)}\right)\right\|^2_{L^2(\Omega)}\nonumber\\&&+C\delta\|\nabla\mathcal{U}_{(n+1)}-\nabla\mathcal{U}_{(n)}\|^2_{L^2(\Omega)}+C\delta^2\|\mathcal{U}_{(n+1)}-\mathcal{U}_{(n)}\|^2_{L^2(\Omega)}\nonumber\\&&+\frac{1}{2}\|\mathcal{U}_{(n+1)}-\mathcal{U}_{(n)}\|^2_{L^2(\Omega)}.
\end{eqnarray}
Applying Lemma \ref{eq:le1} and \ref{eq:le2} to (\ref{eq:h16}), there is
\begin{eqnarray*}
\left(\frac{3}{2}-C\delta^2-C\delta h^{-2}\right)\|\mathcal{U}_{(n+1)}-\mathcal{U}_{(n)}\|^2_{L^2(\Omega)}\leq C\delta h^{11-6p}\|\mathcal{U}_{(n)}-\mathcal{U}_{(n-1)}\|^2_{L^2(\Omega)}.
\end{eqnarray*}
If $\delta$ and $h$ are adequate, we get
\begin{eqnarray}
\label{eq:h9}
\|\mathcal{U}_{(n+1)}-\mathcal{U}_{(n)}\|^2_{L^2(\Omega)}\leq C\|\mathcal{U}_{(n)}-\mathcal{U}_{(n-1)}\|^2_{L^2(\Omega)},\quad \mathrm{with}\;\,0<C<1.
\end{eqnarray}
Iterating the equation (\ref{eq:h9}), we have
\begin{eqnarray*}
\|\mathcal{U}_{(n+1)}-\mathcal{U}_{(n)}\|^2_{L^2(\Omega)}\leq C^n\|\mathcal{U}_{(1)}-\mathcal{U}_{(0)}\|^2_{L^2(\Omega)}\rightarrow 0,n\rightarrow +\infty.
\end{eqnarray*}
Going back to the equation (\ref{eq:h10}), we get
\begin{eqnarray*}
\|\mathcal{Y}_{(n+1)}-\mathcal{Y}_{(n)}\|^2_{L^2(\Omega)}\rightarrow 0,n\rightarrow +\infty,
\end{eqnarray*}
what proves the intended.
\end{proof}

For pratical purpose we can consider a small value $tol$ and use this schemes until
\begin{equation*}
\|\mathcal{Y}_{(n+1)}-\mathcal{Y}_{(n)}\|_{L^2(\Omega)}^2<tol\quad\mathrm{and}\quad\|\mathcal{U}_{(n+1)}-\mathcal{U}_{(n)}\|_{L^2(\Omega)}^2<tol.
\end{equation*}

Let $\{\varphi_1,\cdots,\varphi_n\}$ be the lagrangian base of $\mathcal{S}^h$ associated to the partition $\mathcal{T}_h$. Consider the equations (\ref{eq:ff31}) and (\ref{eq:ff42}), using the representations $u^h(x,t_k)=\sum\limits_{i=1}^{n}U_i(t_k)\varphi_i(x)$, $y^h(x,t_k)=\sum\limits_{i=1}^{n}Y_i(t_k)\varphi_i(x)$ and $w^h(x)=\sum \limits_{j=1}^{n}W_j\varphi_j(x)$, with arbitrary $W_j$, there is
\begin{eqnarray}\label{eq:ff50}
&&\sum\limits_{i=1}^{n}\sum\limits_{j=1}^{n}\int_{\Omega}^{}\left(\frac{U_i(t_{k+1})\varphi_i(x)-U_i(t_{k})\varphi_i(x)}{\delta}\right)W_j\varphi_j(x)dx\nonumber\\&&+\sum\limits_{i=1}^{n}\sum\limits_{j=1}^{n}\int_{\Omega}^{}\left(\left|\frac{U_i(t_{k+1})\varphi_i\prime(x)+U_i(t_{k})\varphi_i\prime(x)}{2}\right|^{p-2}\right.\nonumber\\&&\left.\times\left(\frac{U_i(t_{k+1})\varphi'_i(x)+U_i(t_{k})\varphi_i\prime(x)}{2}\right)W_j\varphi_j\prime(x)\right)dx\nonumber\\&&=\sum\limits_{i=1}^{n}\sum\limits_{j=1}^{n}\int_{\Omega}^{}\left(\frac{Y_i(t_{k+1})\varphi_i(x)+Y_i(t_{k})\varphi_i(x)}{2}\right)W_j\varphi_j(x)dx\nonumber\\&&+\sum\limits_{j=1}^{n}\int_{\Omega}^{}f\big(x,t_{k+\frac{1}{2}}\big)W_j\varphi_j(x)dx.
\end{eqnarray}
Defining the matrices
\begin{eqnarray*}
A(t)\in\mathbb{M}_{n,n}(\mathbb{R}),\quad A(i,j)&=&\int_{\Omega}^{}\left|\sum\limits_{k=1}^{n}U_k(t)\varphi_k\prime(x)\right|^{p-2}\varphi_{i}\prime(x)\varphi_{j}\prime(x)dx, \\M\in\mathbb{M}_{n,n}(\mathbb{R}),\quad M(i,j)&=&\int_{\Omega}^{}\varphi_{i}(x)\varphi_{j}(x)dx,\\F(t)\in\mathbb{M}_{n,1}(\mathbb{R}),\quad F(i)&=&\int_{\Omega}^{}f(x,t)\varphi_{i}(x)dx,
\end{eqnarray*}
with
\begin{equation*}
U(t)=\big(U_1(t),U_2(t),\cdots,U_n(t)\big)^T\quad\mathrm{and}\quad Y(t)=\big(Y_1(t),Y_2(t),\cdots,Y_n(t)\big)^T.
\end{equation*}
Applying matrices to (\ref{eq:ff50}), we get
\begin{eqnarray}\label{eq:ff52}
&&\left(2M+\delta A^{\big(k+\frac{1}{2}\big)}\right)U^{(k+1)}-\delta MY^{(k+1)}\nonumber\\&&=\left(2M-\delta A^{\big(k+\frac{1}{2}\big)}\right)U^{(k)}+\delta MY^{(k)}+2\delta F^{\big(k+\frac{1}{2}\big)}.
\end{eqnarray}

The scheme (\ref{eq:j3}) applied to (\ref{eq:ff52}) is as follows
\begin{eqnarray*}
&&\left(2M+\delta A\big(\mathcal{U}_{(n)}\big)\right)\mathcal{U}_{(n+1)}-\delta M\mathcal{Y}_{(n+1)}\nonumber\\&&=\left(2M-\delta A\big(\mathcal{U}_{(n)}\big)\right)U^{(k)}+\delta MY^{(k)}+2\delta F^{\big(k+\frac{1}{2}\big)}.
\end{eqnarray*}
And the scheme (\ref{eq:h13}) applied to (\ref{eq:ff52}) is
\begin{eqnarray*}
&&2M\mathcal{U}_{(n+1)}-\delta M\mathcal{Y}_{(n+1)}\nonumber\\&&=-\delta A\big(\mathcal{U}_{(n)}\big)\mathcal{U}_{(n)}+\left(2M-\delta A\big(\mathcal{U}_{(n)}\big)\right)U^{(k)}+\delta MY^{(k)}+2\delta F^{\big(k+\frac{1}{2}\big)}.
\end{eqnarray*}
Using the procedures of (\ref{eq:ff50}) in (\ref{eq:ff42}), we have
\begin{eqnarray}\label{eq:ff51}
&&\sum\limits_{i=1}^{n}\sum\limits_{j=1}^{n}\int_{\Omega}^{}\left(\frac{Y_i(t_{k+1})\varphi_i(x)+Y_i(t_{k})\varphi_i(x)}{2}\right)V_j\varphi_j(x)dx\nonumber\\&&=g(0)\sum\limits_{i=1}^{n}\sum\limits_{j=1}^{n}\int_{\Omega}^{}\left(\frac{U_i(t_{k+1})\varphi_i(x)+U_i(t_{k})\varphi_i(x)}{2}\right)V_j\varphi_j(x)dx\nonumber\\&&-g\big(t_{k+\frac{1}{2}}\big)\sum\limits_{i=1}^{n}\sum\limits_{j=1}^{n}\int_{\Omega}^{}U_i(0)\varphi_i(x)V_j\varphi_j(x)dx\nonumber\\&&+\frac{\delta}{2}g\big(t_{k+\frac{1}{2}}\big)\sum\limits_{i=1}^{n}\sum\limits_{j=1}^{n}\int_{\Omega}^{}Y_i(0)\varphi_i(x)V_j\varphi_j(x)dx\nonumber\\&&+\delta\sum\limits_{n=1}^{k-1}g\Big(t_{k+\frac{1}{2}}-t_n\Big)\sum\limits_{i=1}^{n}\sum\limits_{j=1}^{n}\int_{\Omega}^{}Y_i(t_n)\varphi_i(x)V_j\varphi_j(x)dx\nonumber\\&&+\left(\frac{3\delta}{4}g\Big(t_{k+\frac{1}{2}}-t_k\Big)+\frac{\delta}{8}g(0)\right)\sum\limits_{i=1}^{n}\sum\limits_{j=1}^{n}\int_{\Omega}^{}Y_i(t_k)\varphi_i(x)V_j\varphi_j(x)dx\nonumber\\&&+\frac{\delta}{8}g(0)\sum\limits_{i=1}^{n}\sum\limits_{j=1}^{n}\int_{\Omega}^{}Y_i(t_{k+1})\varphi_i(x)V_j\varphi_j(x)dx\nonumber\\&&-\frac{\delta}{2}g\prime\big(t_{k+\frac{1}{2}}\big)\sum\limits_{i=1}^{n}\sum\limits_{j=1}^{n}\int_{\Omega}^{}U_i(0)\varphi_i(x)V_j\varphi_j(x)dx\nonumber\\&&-\delta\sum\limits_{n=1}^{k-1}g\prime\Big(t_{k+\frac{1}{2}}-t_n\Big)\sum\limits_{i=1}^{n}\sum\limits_{j=1}^{n}\int_{\Omega}^{}U_i(t_n)\varphi_i(x)V_j\varphi_j(x)dx\nonumber\\&&-\left(\frac{3\delta}{4}g\prime\Big(t_{k+\frac{1}{2}}-t_k\Big)+\frac{\delta}{8}g\prime(0)\right)\sum\limits_{i=1}^{n}\sum\limits_{j=1}^{n}\int_{\Omega}^{}U_i(t_k)\varphi_i(x)V_j\varphi_j(x)dx\nonumber\\&&-\frac{\delta}{8}g\prime(0)\sum\limits_{i=1}^{n}\sum\limits_{j=1}^{n}\int_{\Omega}^{}U_i(t_{k+1})\varphi_i(x)V_j\varphi_j(x)dx+I(f),
\end{eqnarray}
with
\begin{eqnarray*}
I(f)&\approx &\frac{\delta}{4}g\big(t_{k+\frac{1}{2}}\big)\sum\limits_{j=1}^{n}\int_{\Omega}^{}f(x,t_0)V_j\varphi_j(x)dx\nonumber\\ &&+\frac{3\delta}{4}g\big(t_{k+\frac{1}{2}}-t_{\frac{1}{2}}\big)\sum\limits_{j=1}^{n}\int_{\Omega}^{}f\big(x,t_{\frac{1}{2}}\big)V_j\varphi_j(x)dx\nonumber\\ &&+\delta\sum\limits_{m=1}^{k}g\Big(t_{k+\frac{1}{2}}-t_{m+\frac{1}{2}}\Big)\sum\limits_{j=1}^{n}\int_{\Omega}^{}f\big(x,t_{m+\frac{1}{2}}\big)V_j\varphi_j(x)dx\nonumber\\ &&+\frac{\delta}{2}g(0)\sum\limits_{j=1}^{n}\int_{\Omega}^{}f\big(x,t_{k+\frac{1}{2}}\big)V_j\varphi_j(x)dx.
\end{eqnarray*}
As the constants $V_j$ are arbitrary, applying matrices to (\ref{eq:ff51}), we have
\begin{eqnarray*}
\label{eq:ff53}
&&\left(\frac{1}{2}-\frac{\delta}{8}g(0)\right)MY^{(k+1)}+\left(-\frac{1}{2}g(0)+\frac{\delta}{8}g\prime(0)\right)MU^{(k+1)}\nonumber\\&&=\left(-\frac{1}{2}+\frac{3\delta}{4}g\big(t_{k+\frac{1}{2}}-t_k\big)+\frac{\delta}{8}g(0)\right)MY^{(k)}+\delta\sum\limits_{m=1}^{k-1}g\left(t_{k+\frac{1}{2}}-t_{m}\right)MY^{(m)}\nonumber\\&&+\frac{\delta}{2}g\big(t_{k+\frac{1}{2}}\big)MY^{(0)}+\left(\frac{1}{2}g(0)-\frac{3\delta}{4}g\prime\big(t_{k+\frac{1}{2}}-t_k\big)-\frac{\delta}{8}g(0)\right)MU^{(k)}\nonumber\\&&-\left(g\big(t_{k+\frac{1}{2}}\big)+\frac{\delta}{2}g\prime\big(t_{k+\frac{1}{2}}\big)\right)MU^{(0)}+\delta\sum\limits_{m=1}^{k-1}g\prime\left(t_{k+\frac{1}{2}}-t_{m}\right)MU^{(m)}\nonumber\\&&+\frac{\delta}{4}g\big(t_{k+\frac{1}{2}}\big)F^{(0)}+\frac{3\delta}{4}g\left(t_{k+\frac{1}{2}}-t_{\frac{1}{2}}\right)F^{\big(\frac{1}{2}\big)}+\delta\sum\limits_{m=1}^{k}g\left(t_{k+\frac{1}{2}}-t_{m+\frac{1}{2}}\right)F^{\big(m+\frac{1}{2}\big)}\nonumber\\&&+\frac{\delta}{2}g(0)F^{\big(k+\frac{1}{2}\big)}.
\end{eqnarray*}
The scheme (\ref{eq:j4}) corresponds to
\begin{eqnarray*}
&&\left(-\frac{1}{2}g(0)+\frac{\delta}{8}g\prime(0)\right)M\mathcal{U}_{(n+1)}+\left(\frac{1}{2}-\frac{\delta}{8}g(0)\right)M\mathcal{Y}_{(n+1)}\nonumber\\&&=\left(-\frac{1}{2}+\frac{3\delta}{4}g\big(t_{k+\frac{1}{2}}-t_k\big)+\frac{\delta}{8}g(0)\right)MY^{(k)}+\delta\sum\limits_{m=1}^{k-1}g\left(t_{k+\frac{1}{2}}-t_m\right)MY^{(m)}\nonumber\\&&+\frac{\delta}{2}g\big(t_{k+\frac{1}{2}}\big)MY^{(0)}+\left(\frac{1}{2}g(0)-\frac{3\delta}{4}g\prime\big(t_{k+\frac{1}{2}}-t_k\big)-\frac{\delta}{8}g(0)\right)MU^{(k)}\nonumber\\&&-\left(g\big(t_{k+\frac{1}{2}}\big)+\frac{\delta}{2}g\prime\big(t_{k+\frac{1}{2}}\big)\right)MU^{(0)}+\delta\sum\limits_{m=1}^{k-1}g\prime\left(t_{k+\frac{1}{2}}-t_m\right)MU^{(m)}\nonumber\\&&+\frac{\delta}{4}g\big(t_{k+\frac{1}{2}}\big)F^{(0)}+\frac{3\delta}{4}g\left(t_{k+\frac{1}{2}}-t_{\frac{1}{2}}\right)F^{\big(\frac{1}{2}\big)}+\delta\sum\limits_{m=1}^{k}g\left(t_{k+\frac{1}{2}}-t_{m+\frac{1}{2}}\right)F^{\big(m+\frac{1}{2}\big)}\nonumber\\&&+\frac{\delta}{2}g(0)F^{\big(k+\frac{1}{2}\big)}.
\end{eqnarray*}

\section{Numerical results}

In this section we make several simulations in order to study the error, stability and convergence of the method and also to analyze the influence of the parameter and the memory term in the properties of the solutions.

\subsection*{Example 1}

In this example, the study of the convergence of $h$ and $\delta$ using Lagrange basis of degree $1,\,2,\,3$ and $4$ is presented. Since the explicit solutions provided by Antonsev et al. are not easy to calculate, we define an appropriate solution and calculate the corresponding $f$ function. We considered $\Omega=\,]0,1[$, $\lambda=1$, $T=0.1$, $u(x,t)=\big(x(1-x)\big)^2e^{-t}$, $g(\xi)=\lambda e^{-\xi}$, the tolerance of the fixed point method is $tol=10^{-9}$ and $f$ is calculated such that $u$ is the exact solution.

\subsubsection*{Case 1}

\begin{center}
\begin{figure}[h]
\centering
\includegraphics[height=0.182\paperheight]{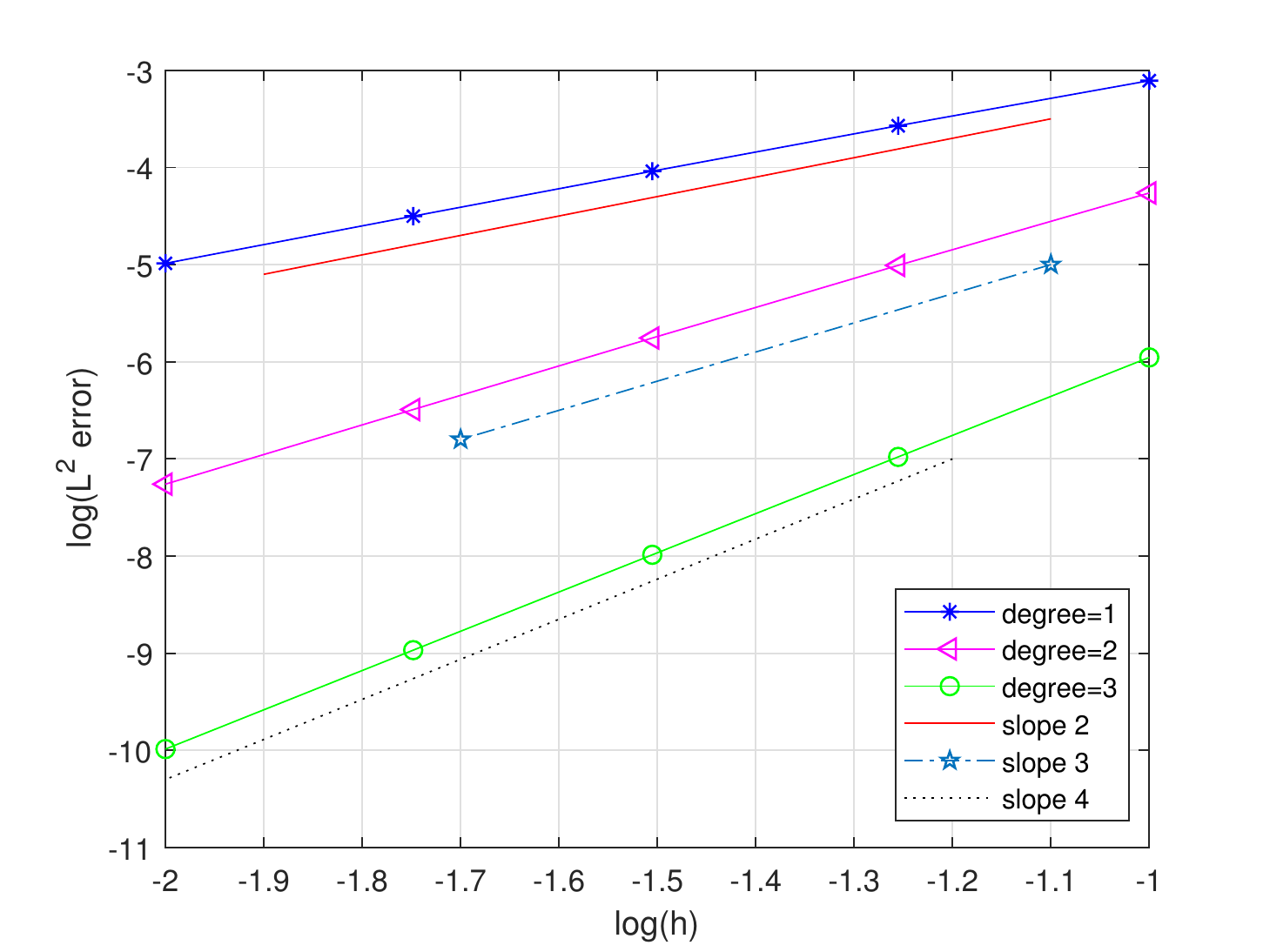}\hfill
\includegraphics[height=0.182\paperheight]{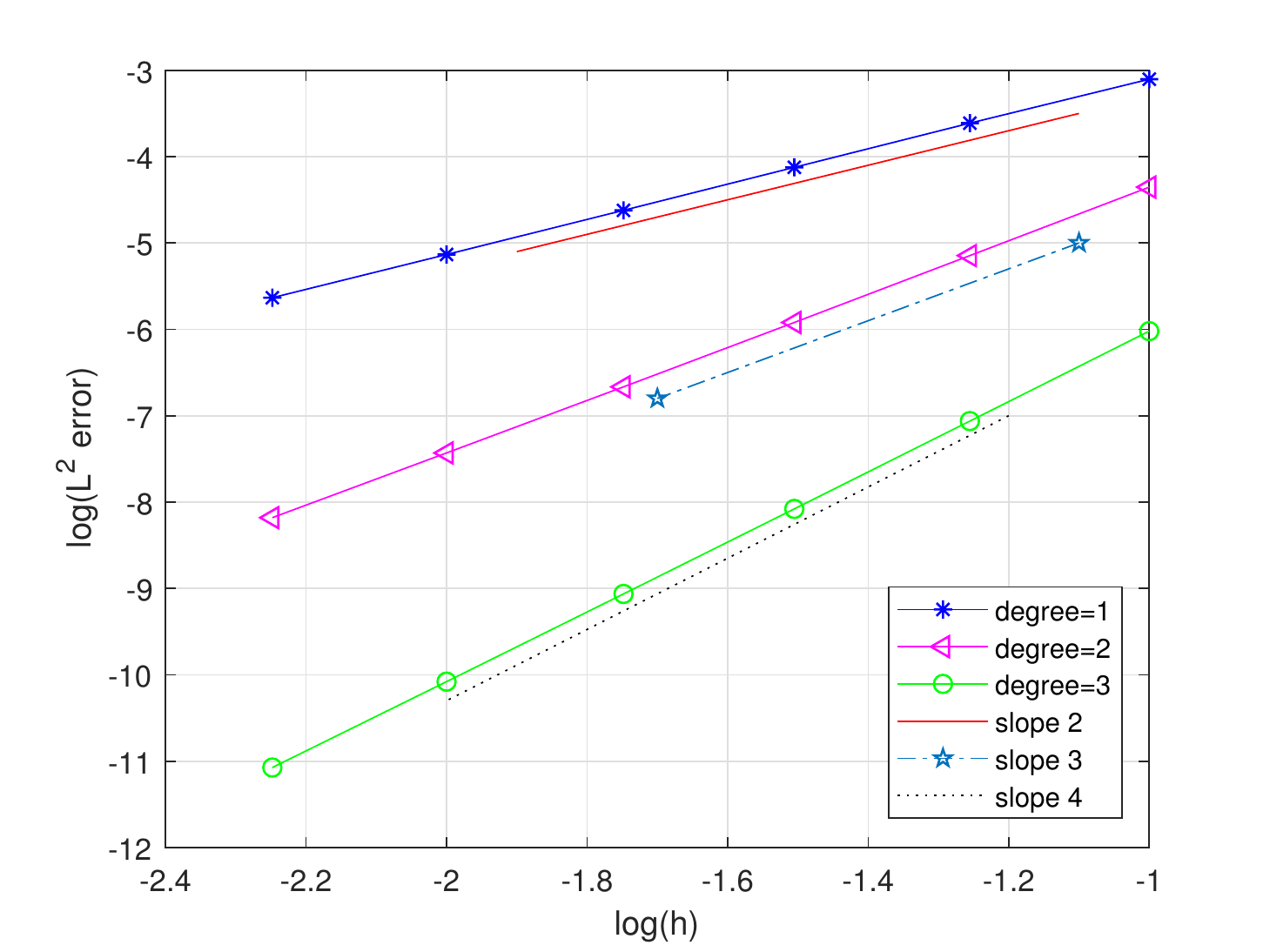}
\caption{Convergence order for $r=1,\,2,\,3$ for $p=3$ (left) and $p=4$ (right). The convergence observed is in agreement with Theorem \ref{the4}.}
\label{f1}
\end{figure}
\end{center}

First we study the convergence for $h$. Fixing a small $\delta =10^{-5}$, varying the values of $h$ and calculating the corresponding $L^2$ error, with the formula $L^2\,error=\left(\int_{\Omega}^{}\big(u(x,T)-U^N(x)\big)^2\right)^{\frac{1}{2}}$ we build the graphs in figure \ref{f1}. We can observe that the method converges and the numerical order of convergence is $2$ for polynomials of degree $1$, $3$ for degree $2$ and $4$ for degree $3$.

\subsubsection*{Case 2}

\begin{center}
\begin{figure}[h]
\centering
\includegraphics[height=0.20\paperheight]{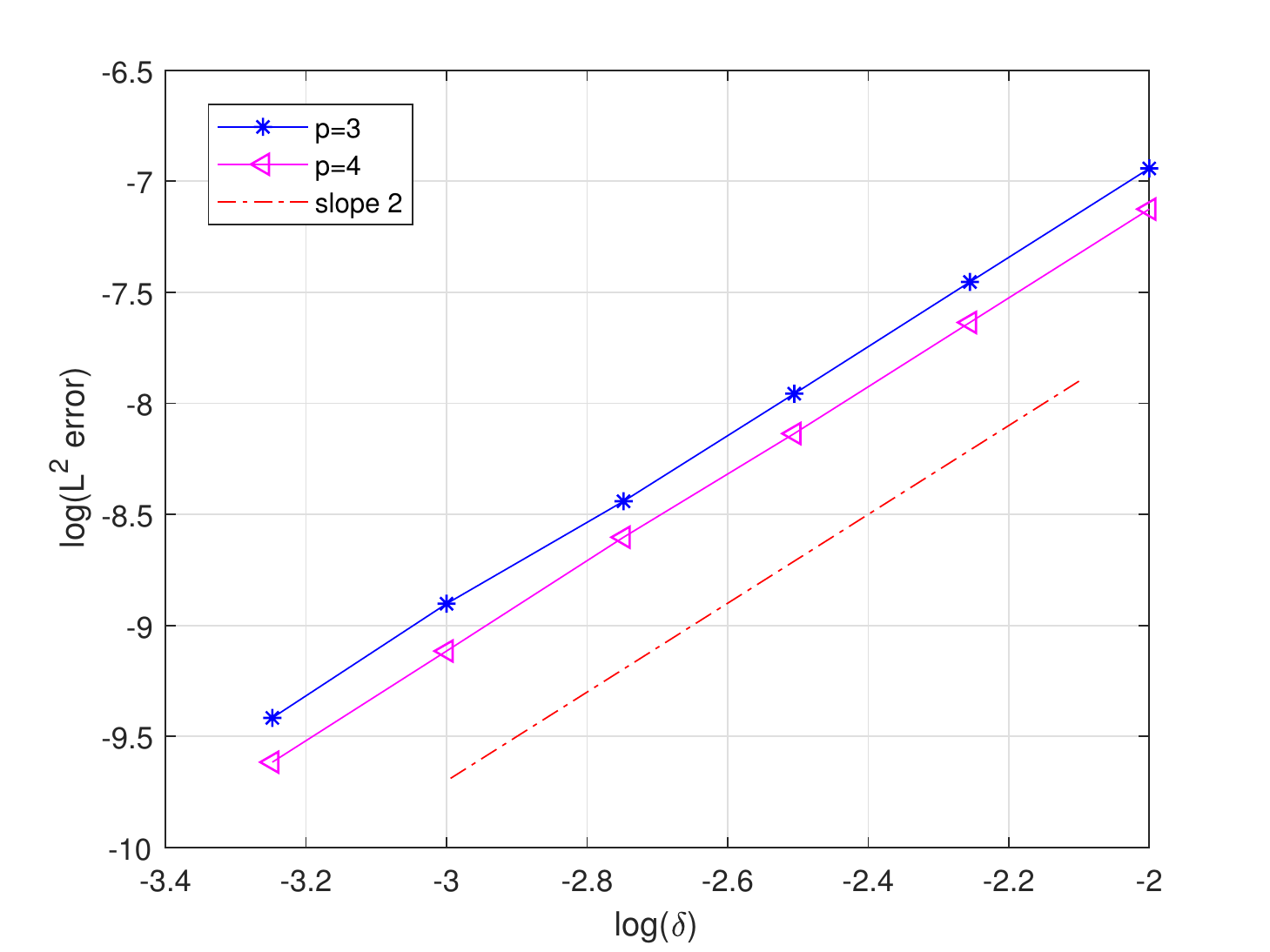}
\caption{Convergence order $\delta$ for $p=3$ and $p=4$.}
\label{f2}
\end{figure}
\end{center}

For the fixed-point method to converge $h$ and $\delta$ must satisfy a certain condition, see Theorem \ref{the1} and \ref{the2}, and indeed for some combinations of $h$ and $\delta$ the method did not converged in the example. Therefore, we cannot fix a small $h$ and vary $\delta$. However, since the solution is a polynomial of degree $4$, if we define the degree of the approximate solution also $4$, the solution will be exact for $x$ and we can measure the error for $\delta$. Fixing $r=4$ and $h=10^{-1}$ and varying $\delta$ we constructed the graph in figure \ref{f2}. Thus, from the graph in figure \ref{f2} we conclude that the order of convergence is $2$, for $p=3$ and $4$ as expected.
\newpage

\subsection*{Example 2}

In this example, we study the asymptotic behavior of solution $u(x)$ for several cases, as described below. Here we consider $u_0(x)=1-x^4$, $g(\xi)=\lambda e^{-\xi}$, $f=0$, $\Omega=[-1,1]$, $r=1$, $T=3$, $tol=10^{-9}$, $h=0.2$ and $\delta=0.001$. To better understand the behaviour we define the energy function $b(t)=\int_{-1}^{1}U(x,t)^2dx$. We simulate with $p=1.5$, $2$ and $4$.

\subsubsection*{Case 1}

\begin{center}
\begin{figure}[h]
\centering
\includegraphics[height=0.182\paperheight]{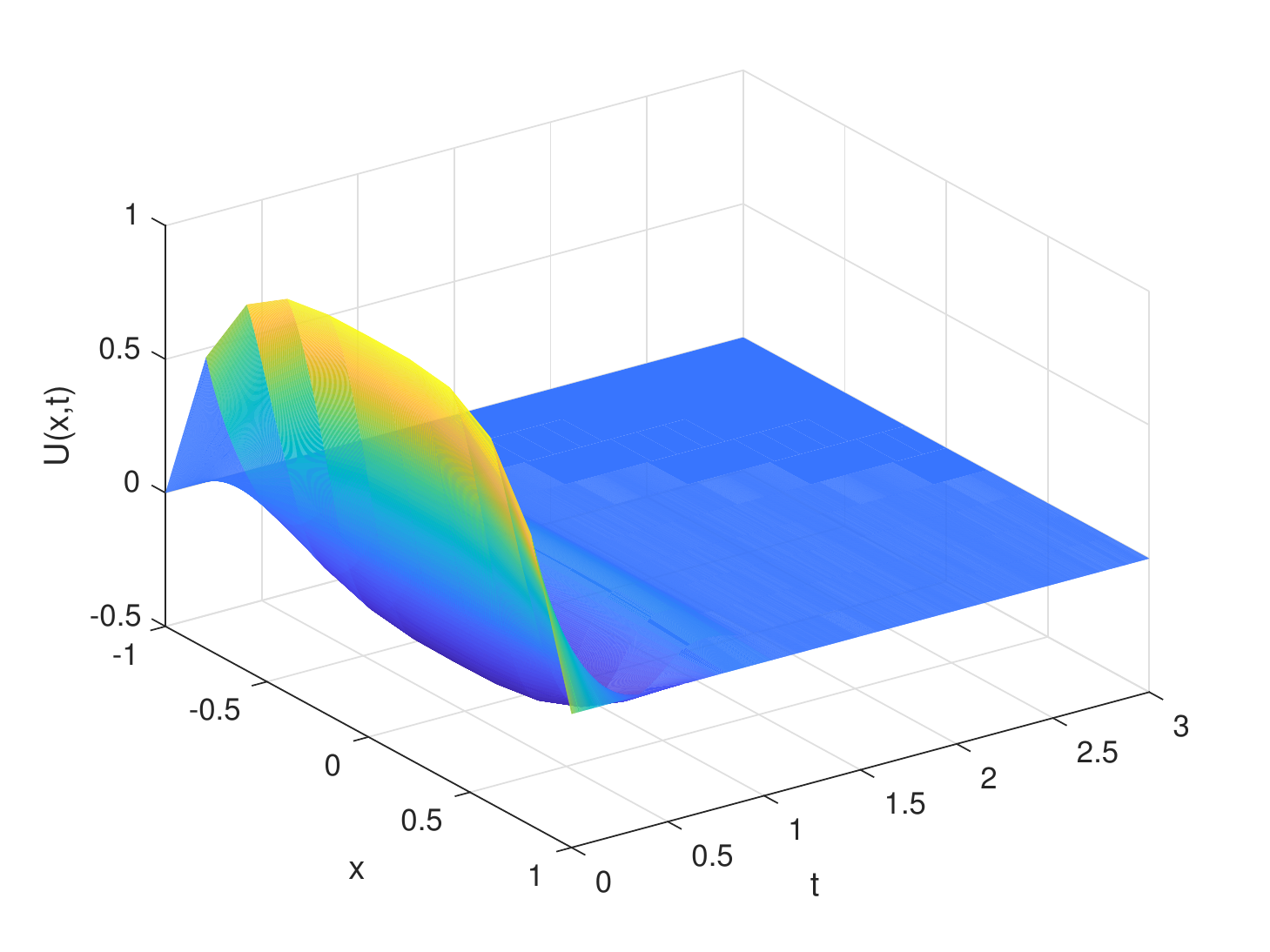}\hfill
\includegraphics[height=0.182\paperheight]{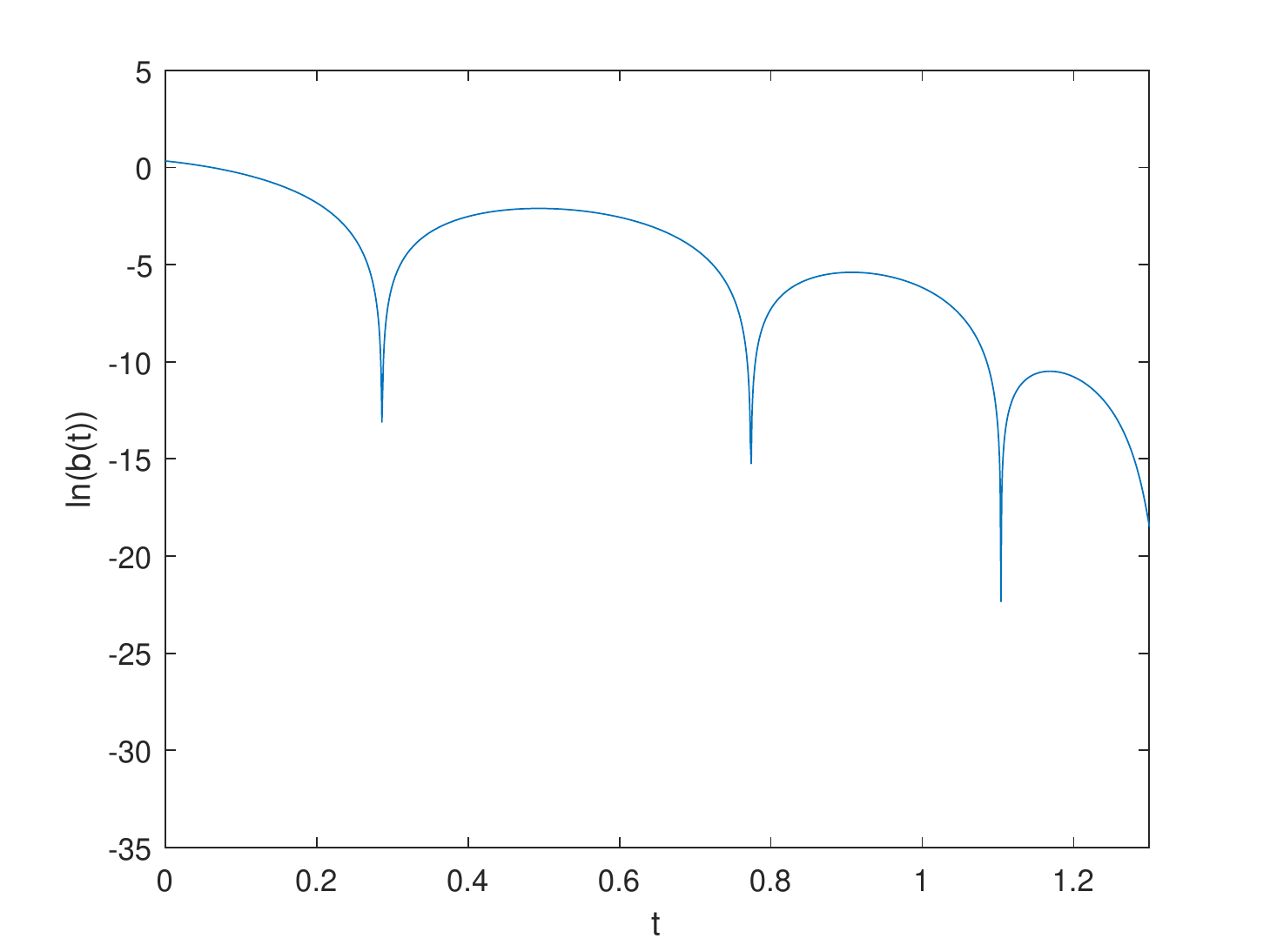}
\caption{Solution obtained (left) and the function $b(t)$ (right) with $p=1.5$.}
\label{f3}
\end{figure}
\end{center}

\begin{center}
\begin{figure}[h]
\centering
\includegraphics[height=0.182\paperheight]{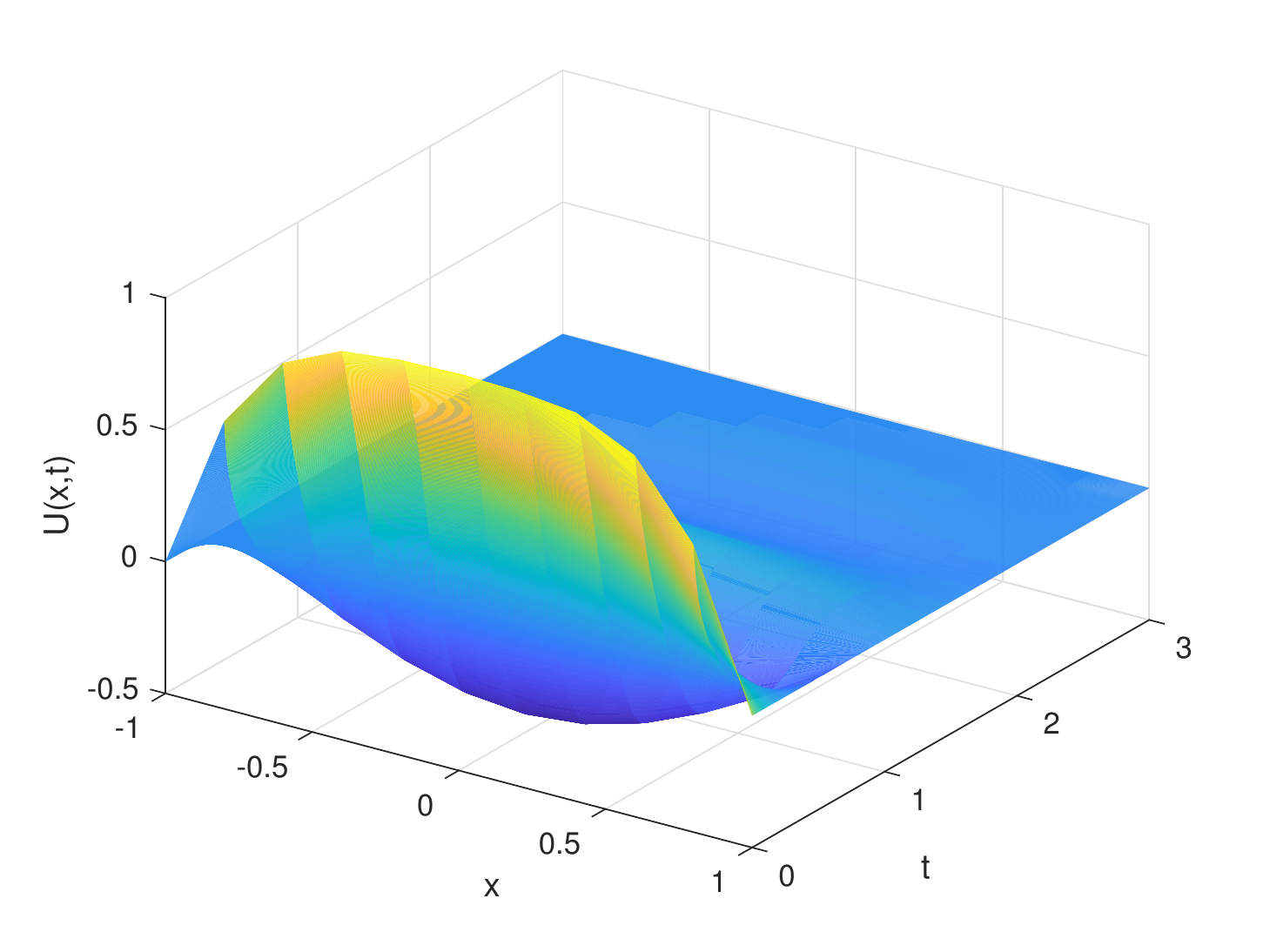}\hfill
\includegraphics[height=0.182\paperheight]{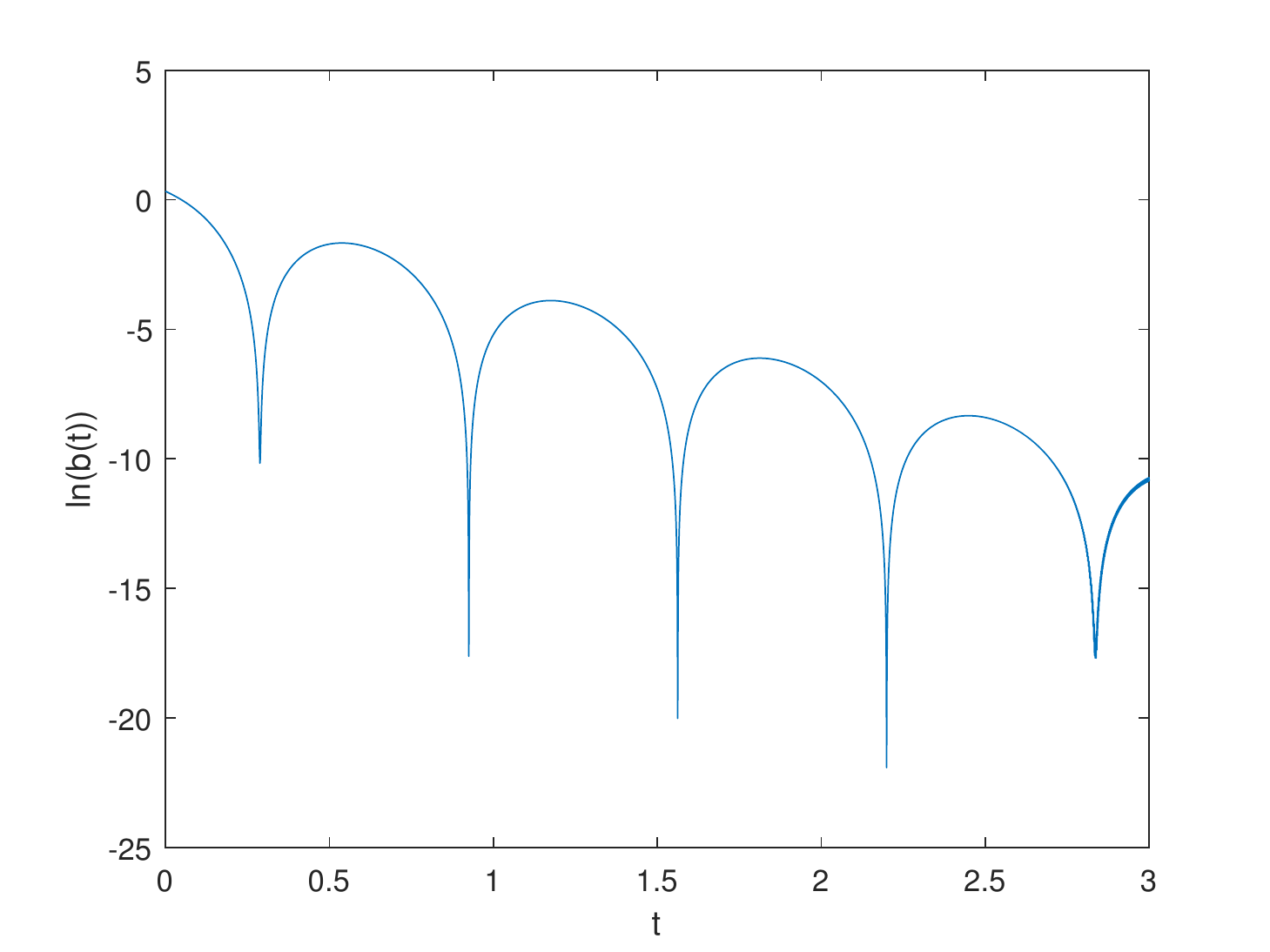}
\caption{Solution obtained (left) and the function $b(t)$ (right) with $p=2$.}
\label{f4}
\end{figure}
\end{center}

\begin{center}
\begin{figure}[h]
\centering
\includegraphics[height=0.182\paperheight]{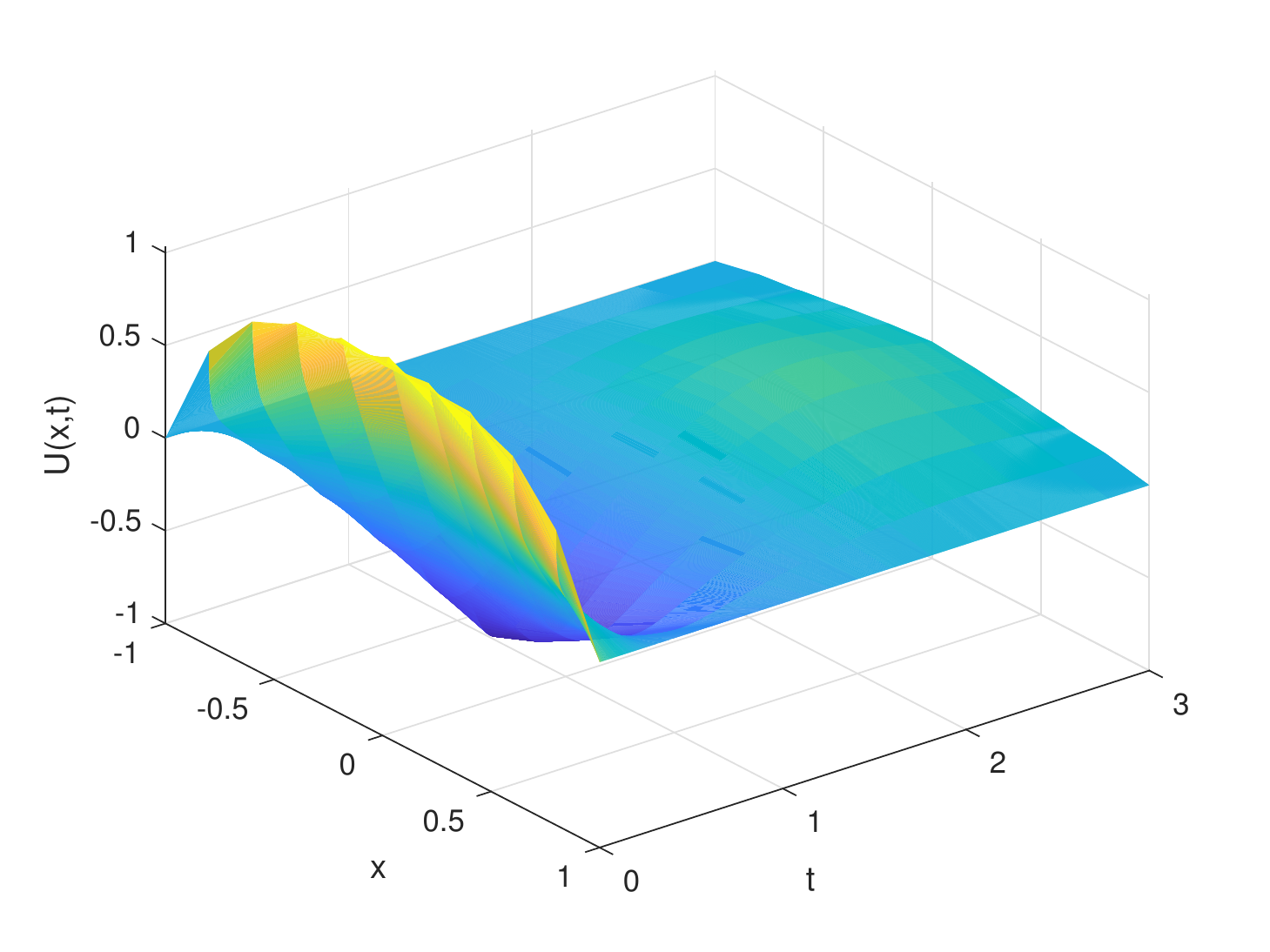}\hfill
\includegraphics[height=0.182\paperheight]{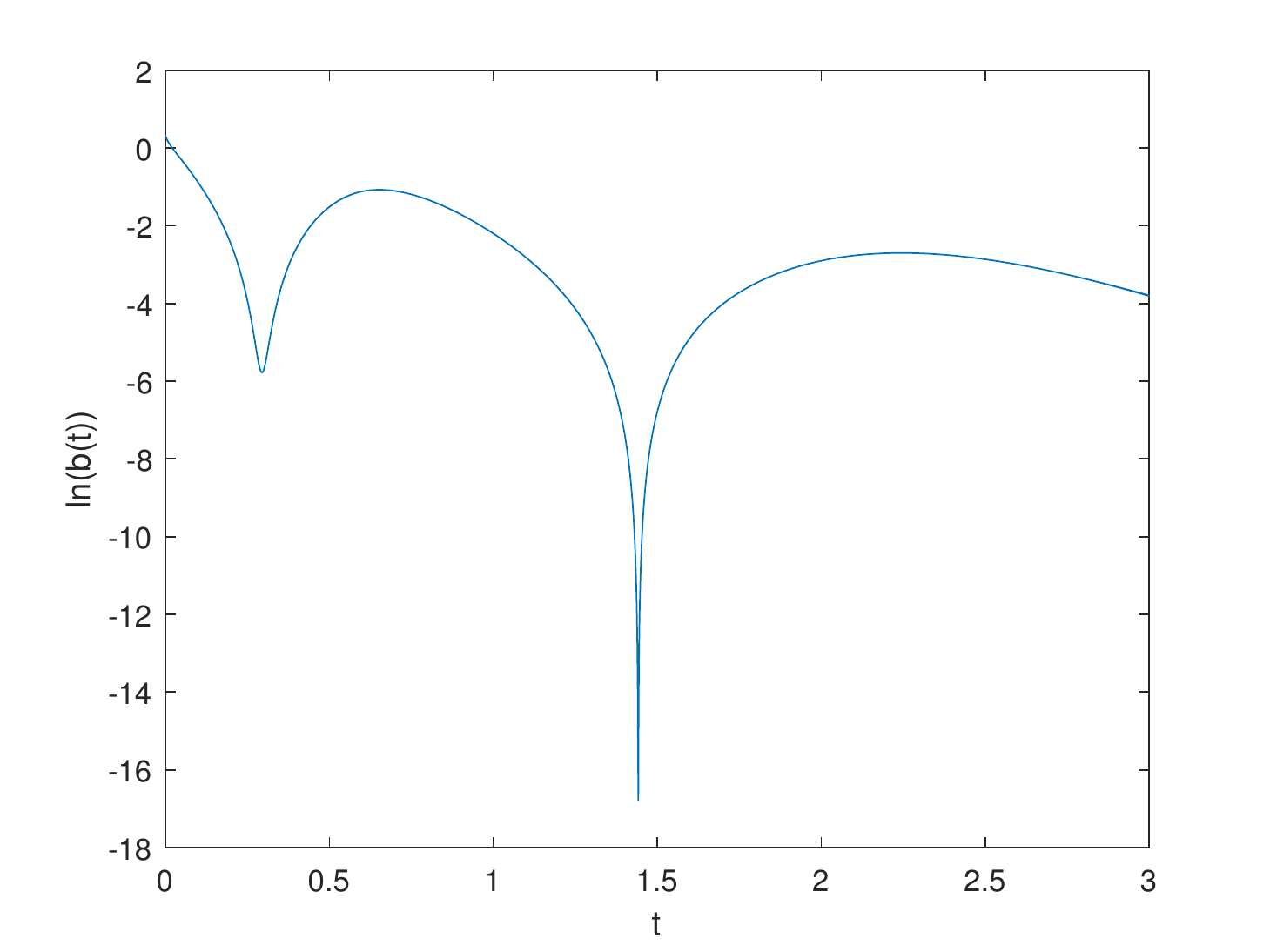}
\caption{Solution obtained (left) and the function $b(t)$ (right) with $p=4$.}
\label{f5}
\end{figure}
\end{center}

For positive $\lambda$ ($\lambda=10$ in this case), the solution decays and develops oscillations in the time direction violating the maximum principle, as we can see in figures \ref{f3}, \ref{f4} and \ref{f5}. The speed of decay depends on the value of $p$.
\newpage

\newpage

\subsubsection*{Case 2}

\begin{center}
\begin{figure}[h]
\centering
\includegraphics[height=0.182\paperheight]{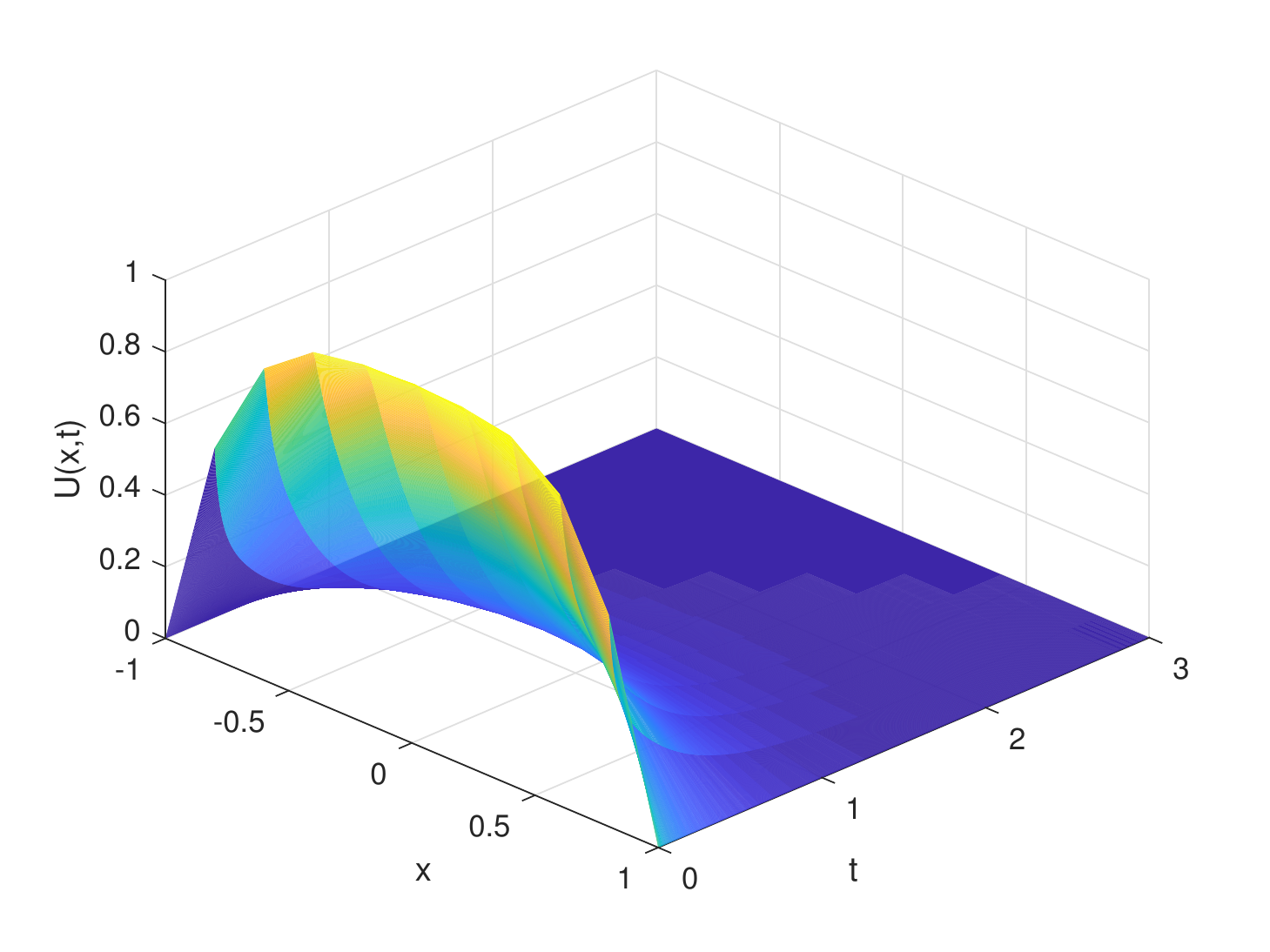}\hfill
\includegraphics[height=0.182\paperheight]{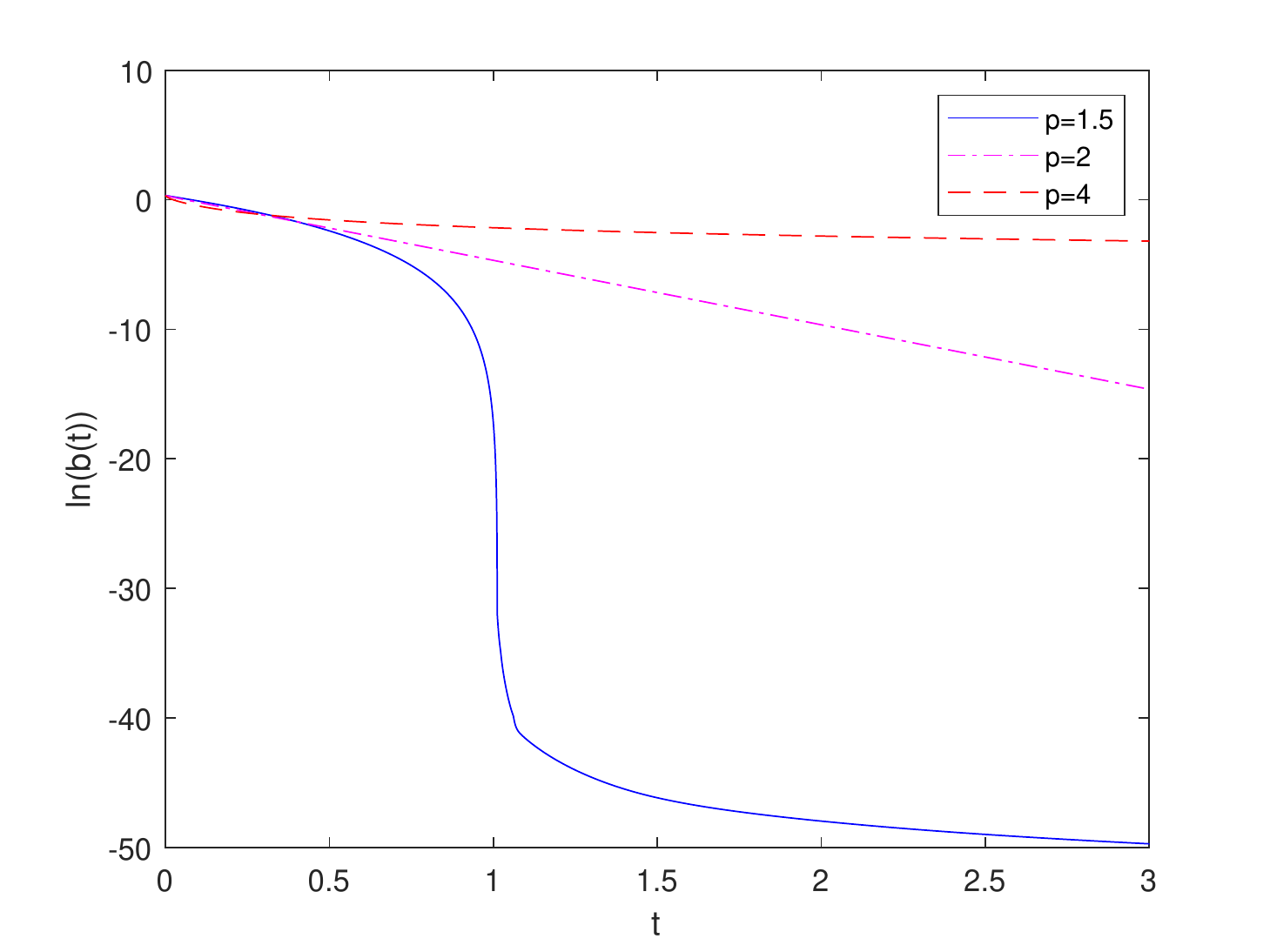}
\caption{Solution obtained (left) and the function $b(t)$ (right) for case $2$.}
\label{f6}
\end{figure}
\end{center}

For the heat equation $(p=2)$ and absence of memory $(\lambda=0)$, the exponential decay is well known, for the case $p=1.5$ we observed an extinction, while for $p=4$ we had a decay, as illustrated in figure \ref{f6}.

\subsubsection*{Case 3}

\begin{center}
\begin{figure}[h]
\centering
\includegraphics[height=0.182\paperheight]{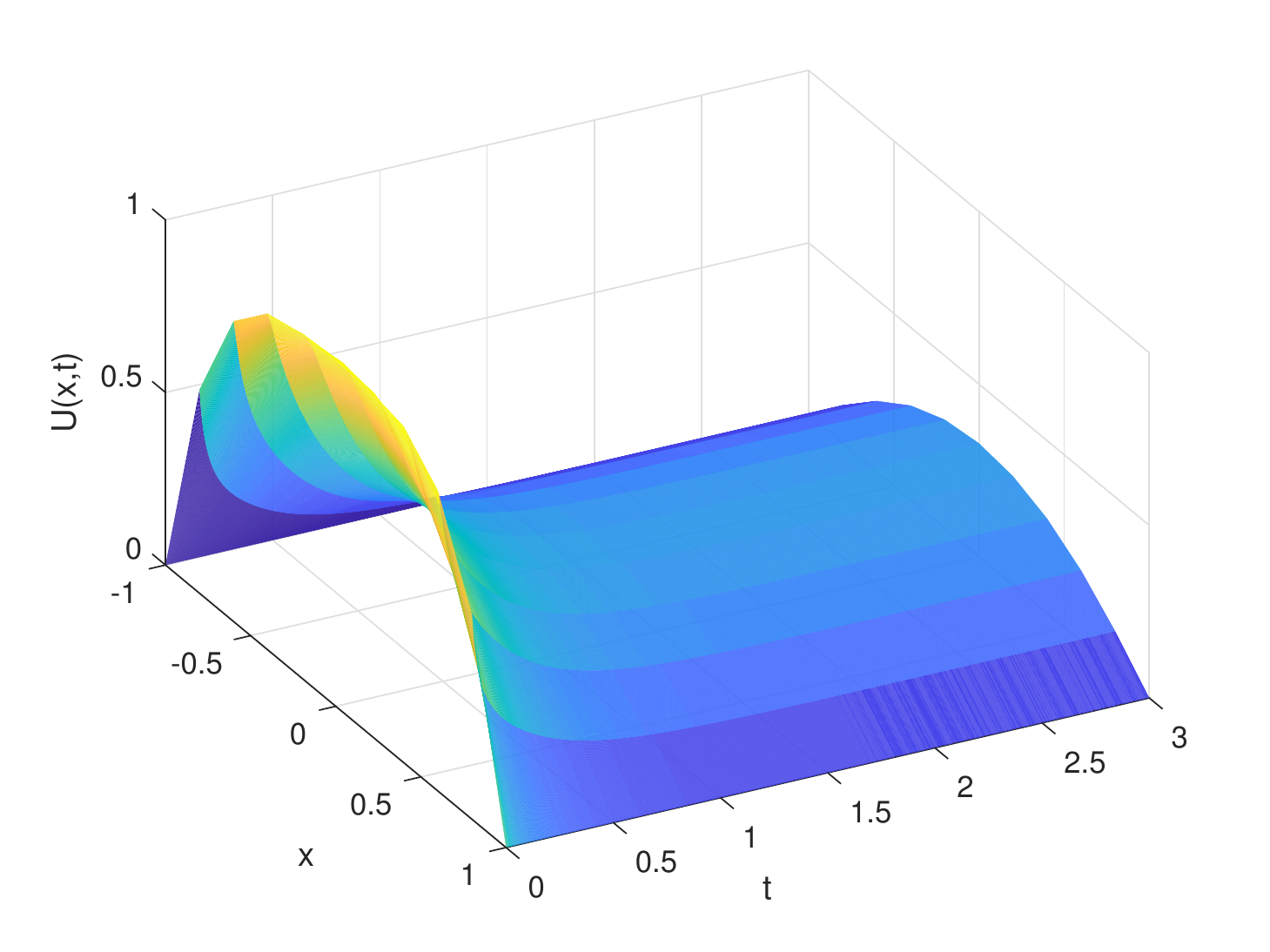}\hfill
\includegraphics[height=0.182\paperheight]{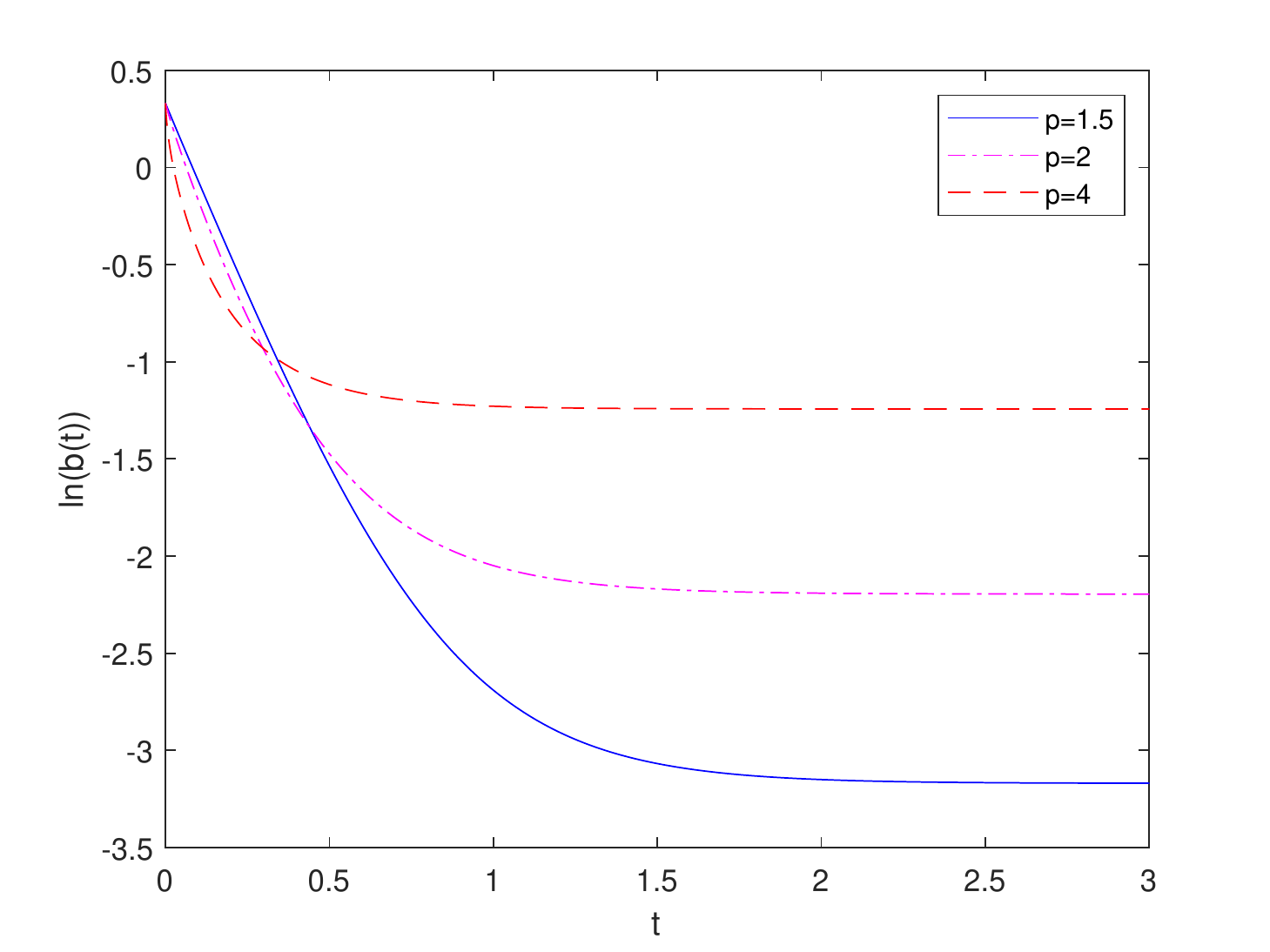}
\caption{Solution obtained (left) and the function $b(t)$ (right) for case $3$.}
\label{f7}
\end{figure}
\end{center}

For negative $\lambda$ with small modulus (in this case $\lambda=-1$) the solution decays with time and approaches a stable solution. However, the speed of decay and the asymptote stable depends on $p$, as illustrated in figure \ref{f7}.

\subsubsection*{Case 4}

\begin{center}
\begin{figure}[h]
\centering
\includegraphics[height=0.182\paperheight]{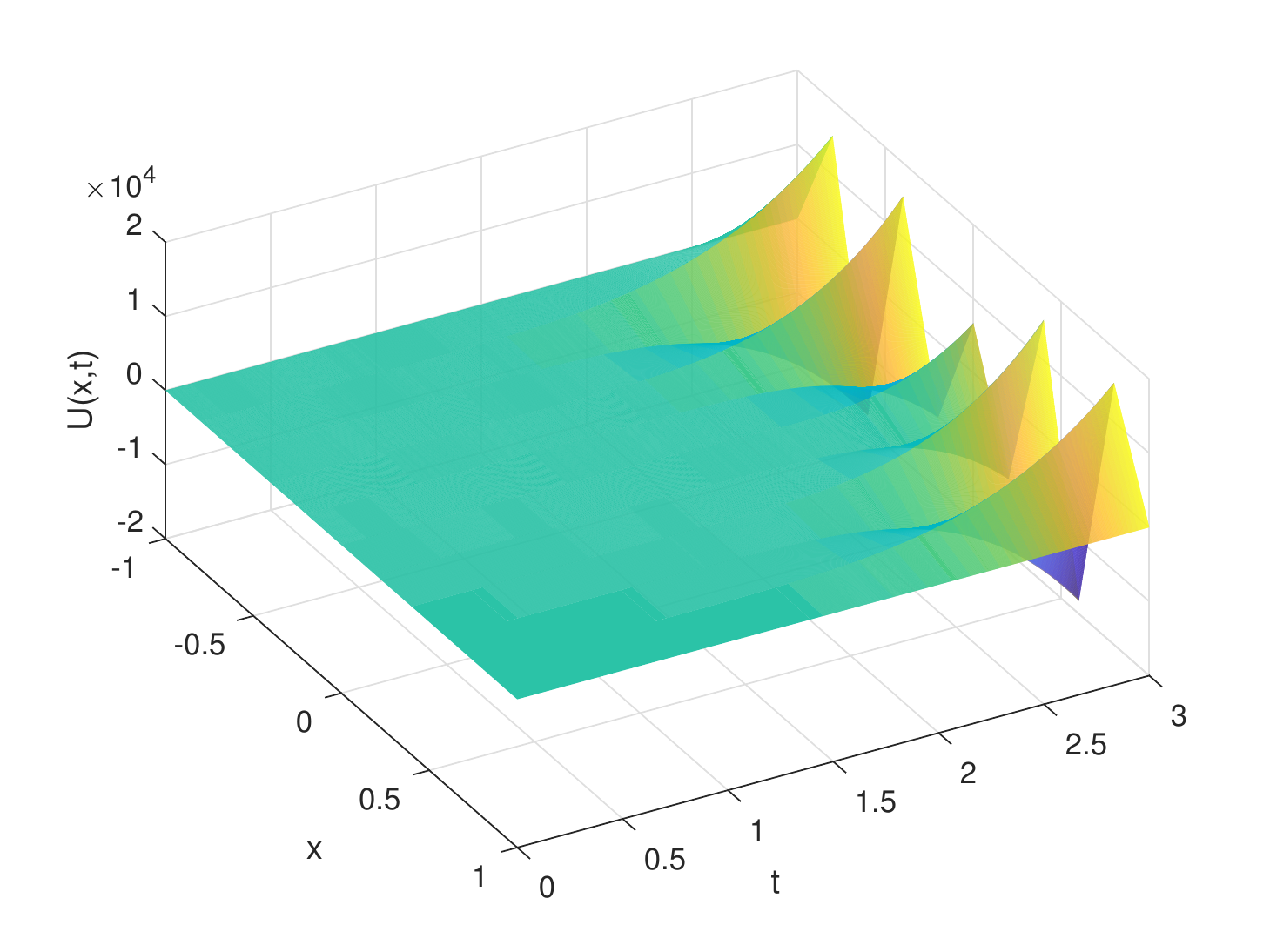}\hfill
\includegraphics[height=0.182\paperheight]{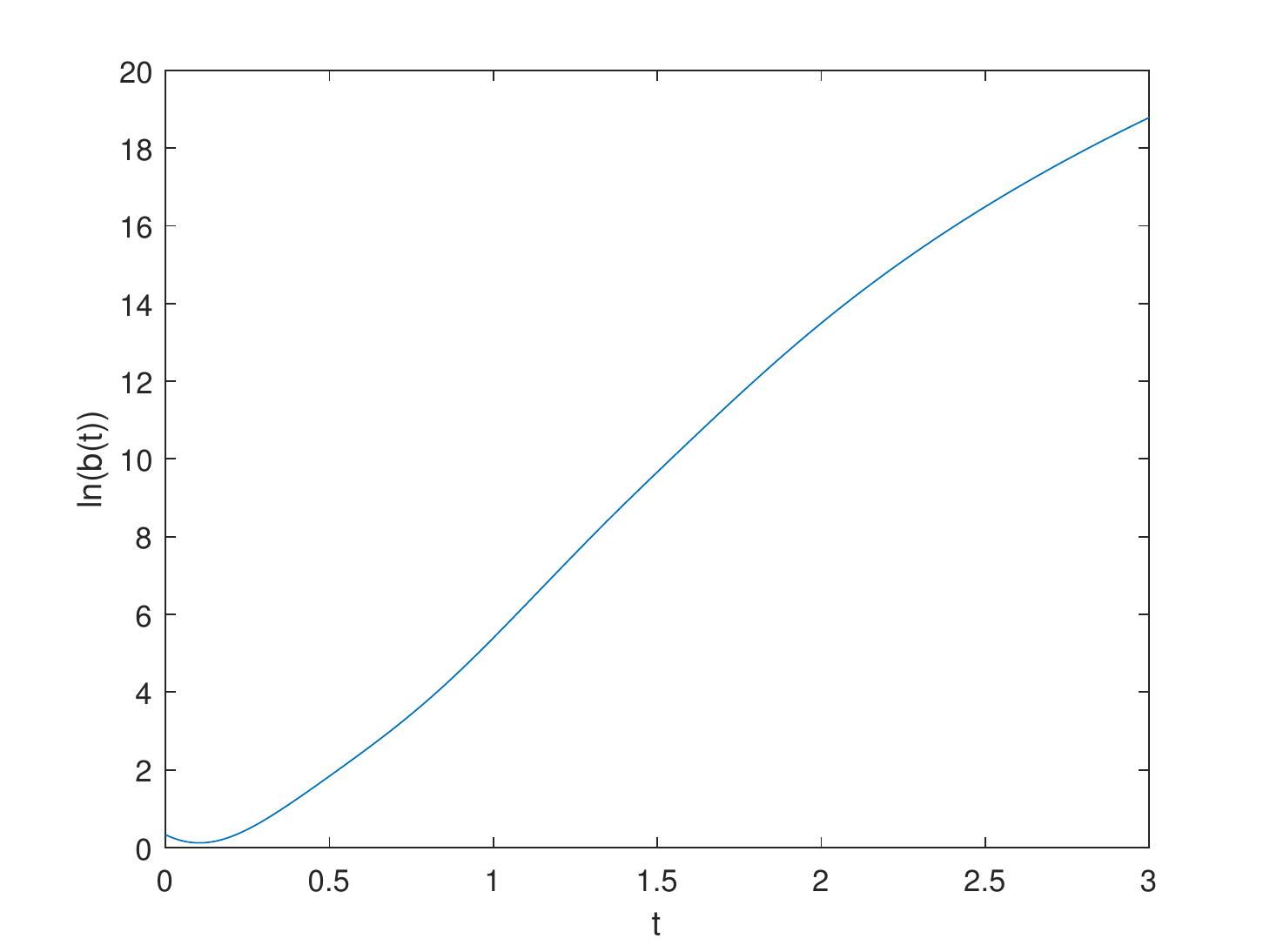}
\caption{Solution obtained (left) and the function $b(t)$ (right) with $p=1.5$ case $4$.}
\label{m1}
\end{figure}
\end{center}

\begin{center}
\begin{figure}[h]
\centering
\includegraphics[height=0.182\paperheight]{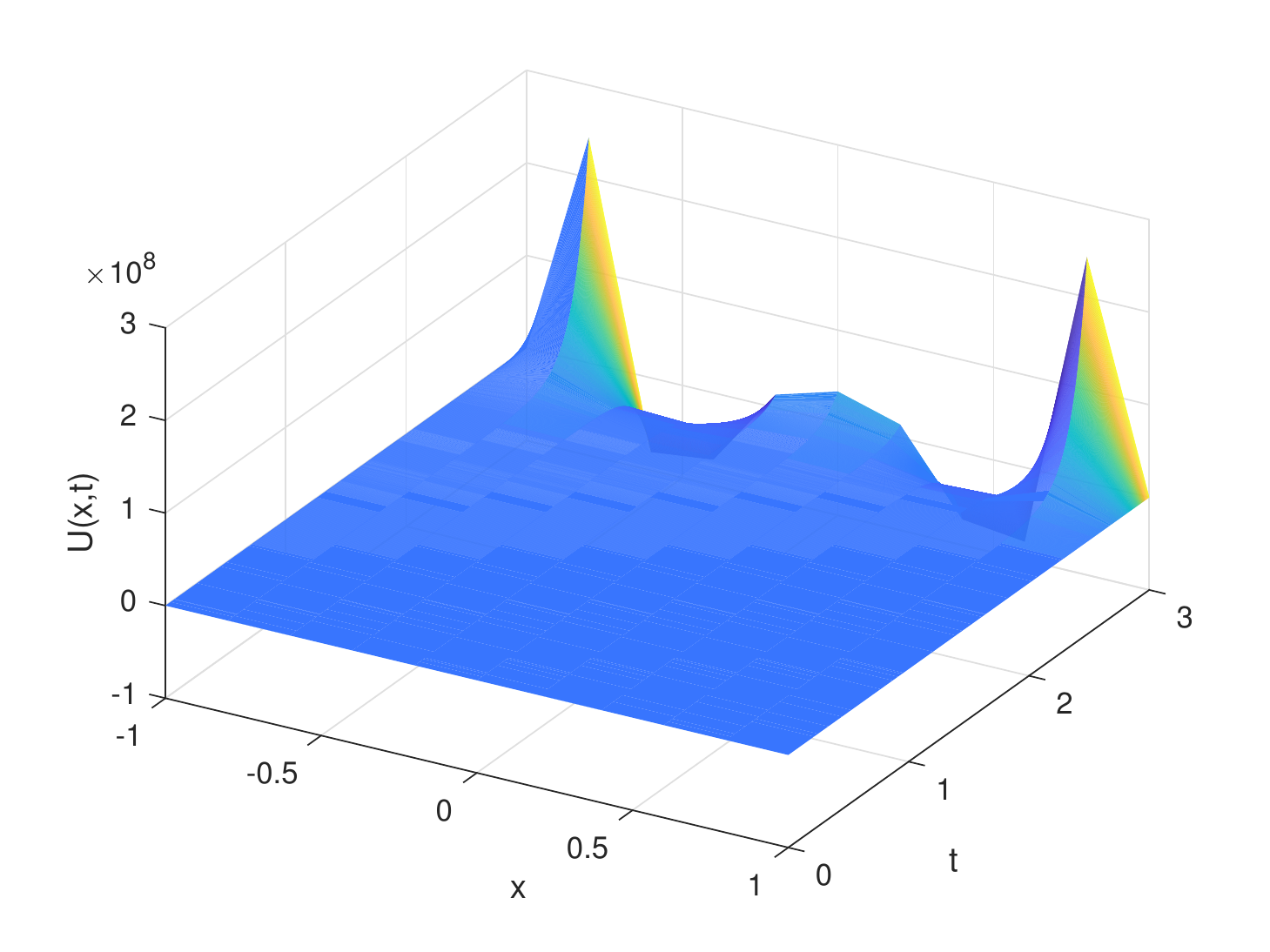}\hfill
\includegraphics[height=0.182\paperheight]{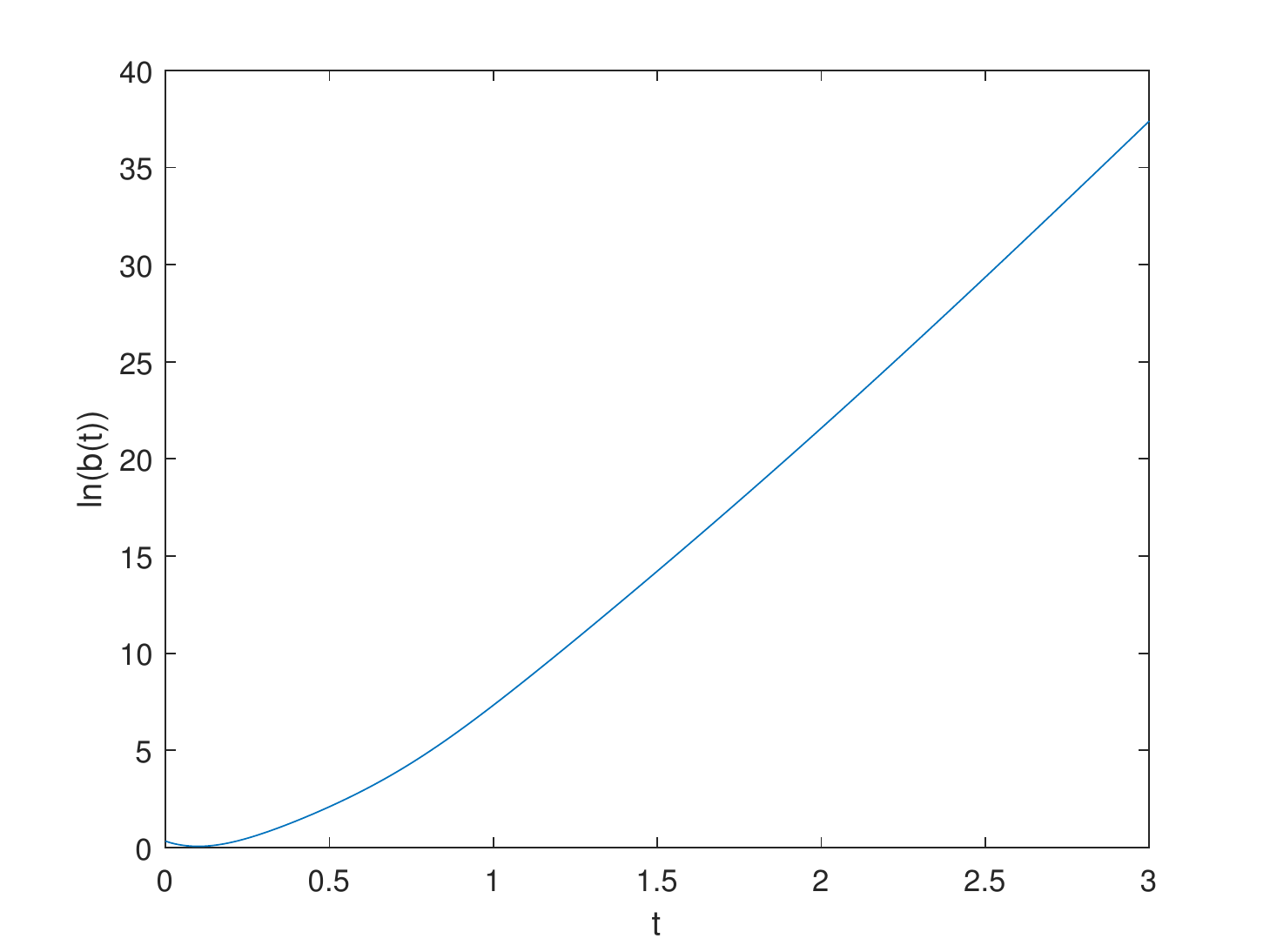}
\caption{Solution obtained (left) and the function $b(t)$ (right) with $p=2$ case $4$.}
\label{m2}
\end{figure}
\end{center}

For negative $\lambda$ with large modulus (in this case $\lambda=-10$) the solution develops oscillations in the direction of space that increase with time, violating also the maximum principle, as we can see in figures \ref{m1} and \ref{m2}.
\newpage

\subsection*{Example 3}

\begin{center}
\begin{figure}[h]
\centering
\includegraphics[height=0.160\paperheight]{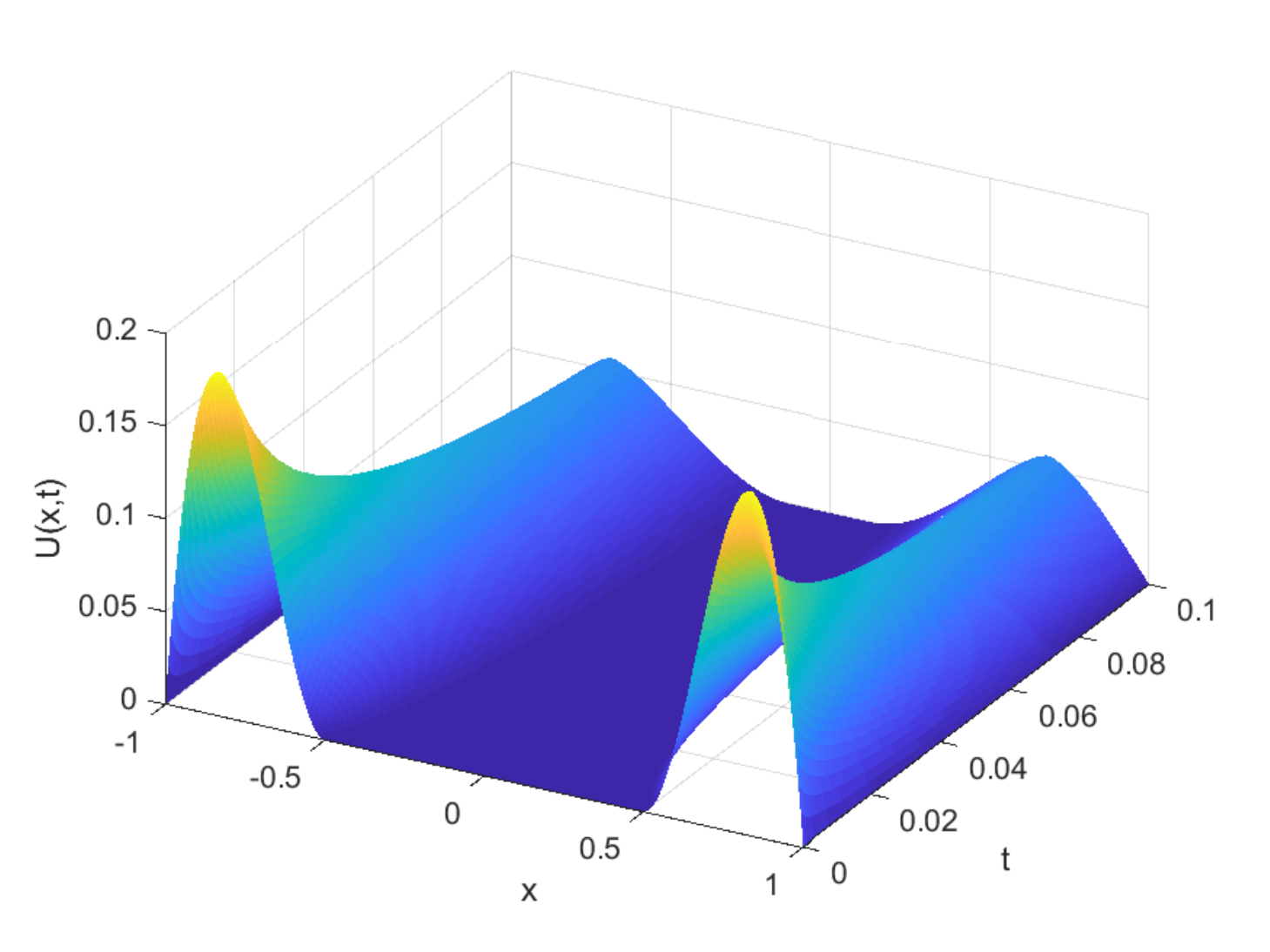}\hfill
\includegraphics[height=0.154\paperheight]{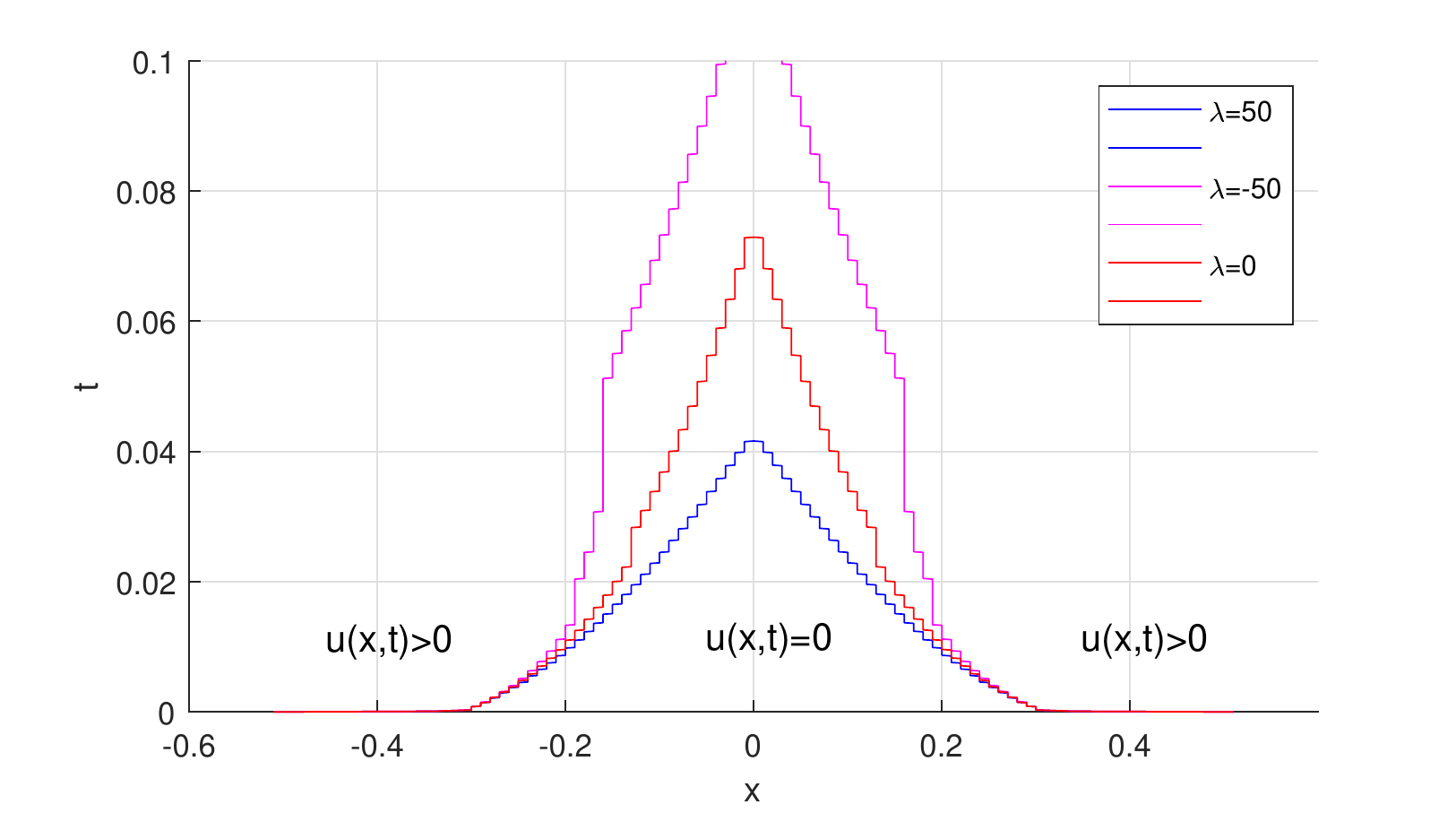}
\caption{Solution obtained (left) and boundary evolution (right) in example $3$.}
\label{f8}
\end{figure}
\end{center}

In this example, we study the finite speed of propagation of solution $u(x)$ for several cases. In \cite{MR3462616} it is proved that for $p>2$, $g\in L^p(0,T)$ and $f=0$, then there is a $t^*>0$ and a positive function $\rho(t)$ such that if $u(x,0)=0$ in $B_{\rho(0)}(x_0)$ then $u(x,t)=0$ in $B_{\rho(t)}(x_0)$ $\forall t\in]0,t^*[$. The data is $g(\xi)=\lambda e^{-\xi}$, $tol=10^{-9}$, $h=0.02$, $\delta=0.001$, $p=3$, $r=1$ and
\begin{eqnarray*}
u_0(x)=
\begin{cases}
10(x+1)(0.5+x)^2, & x\in [-1,-0.5[,\\
0, & x\in [0.5, -0.5],\\
10(1-x)(x-0.5)^2, & x\in ]0.5,1].
\end{cases}
\end{eqnarray*}
We know perfectly well that the heat equation does not exhibit finite speed os propagation. However, for $p=3$ with no memory that is $\lambda=0$, we observe this effect. Starting from this initial data, which is zero between $-0.5$ and $0.5$, it is verified in figure \ref{f8} that the size of the region where the solution is zero is decreasing with finite speed and this speed depends on $\lambda$.

\subsection*{Example 4}

\begin{center}
\begin{figure}[h]
\centering
\includegraphics[height=0.160\paperheight]{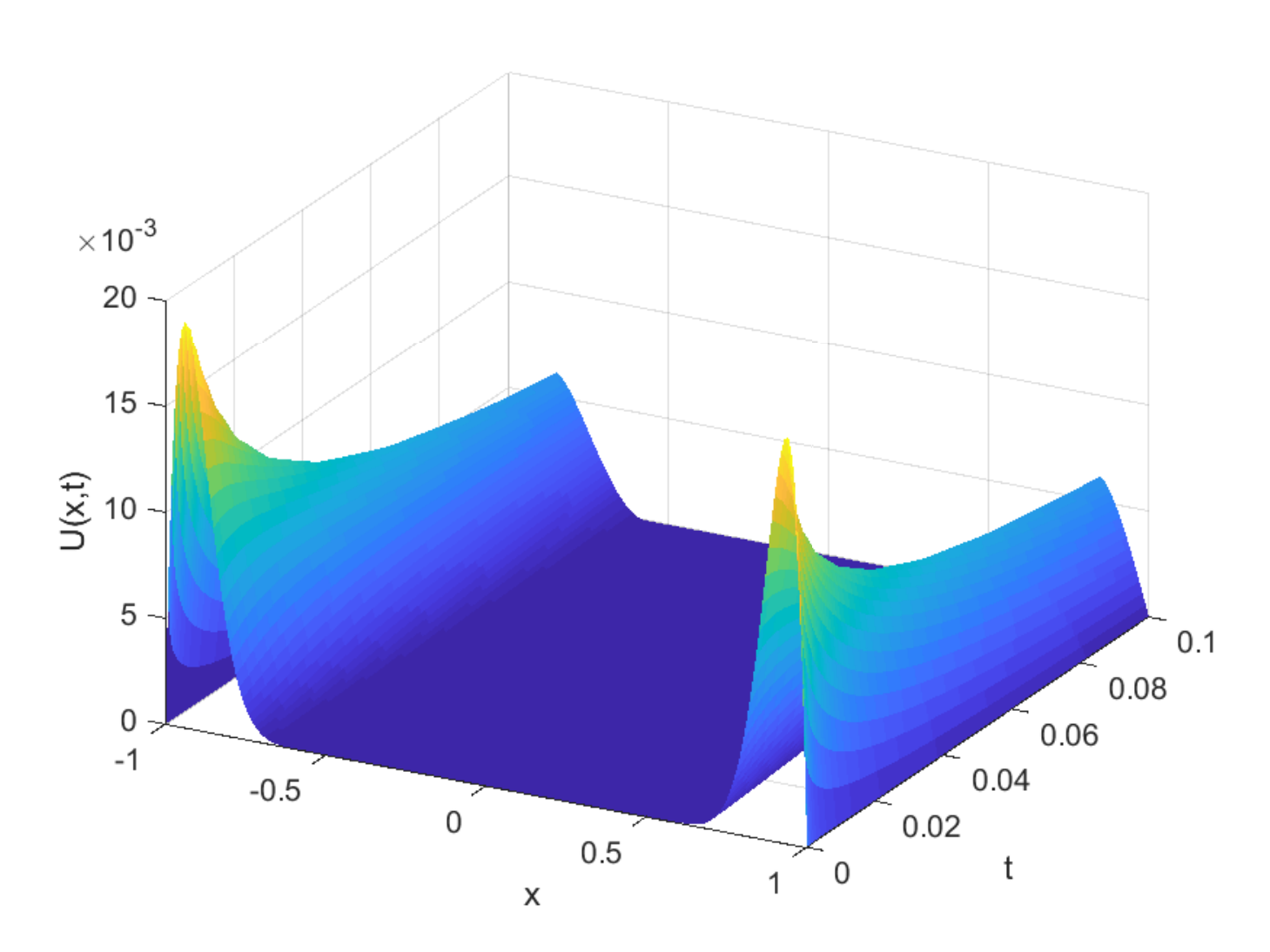}\hfill
\includegraphics[height=0.154\paperheight]{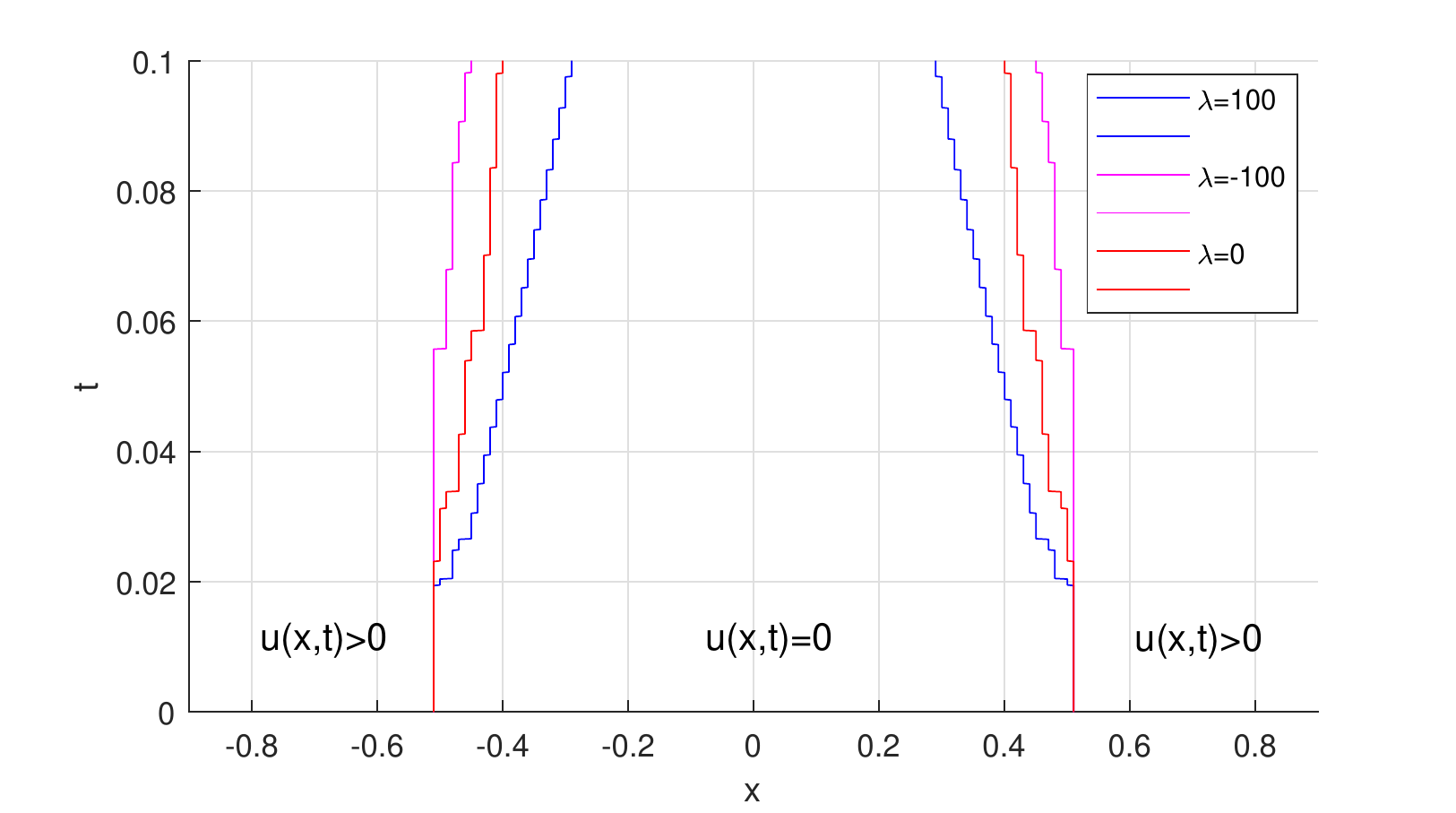}
\caption{Solution obtained (left) and boundary evolution (right) in example $4$.}
\label{f9}
\end{figure}
\end{center}

Finally changing the initial solution around $0.5$ we could observe the waiting time effect. In \cite{MR3462616} it is proved that for $p>2$, if $u_0$ is sufficiently flat near the boundaries of his support then there is $t_*\in\,]0,T[$ such that $u(x,t)=0$ in $B_{\rho_0}\times]0,t^*[$. Let us consider $g(\xi)=\lambda e^{-\xi}$, $tol=10^{-9}$, $h=0.02$, $\delta=0.001$, $p=3$, $r=1$ and
\begin{eqnarray*}
u_0(x)=
\begin{cases}
100(x+1)(0.5+x)^7, & x\in [-1,-0.5[,\\
0, & x\in [0.5, -0.5],\\
100(1-x)(x-0.5)^7, & x\in ]0.5,1].
\end{cases}
\end{eqnarray*}
Comparing with the previous example we observe that the size of the region where the solution is zero remains fixed for a short period of time and after that it decreases with finite speed, as illustrated in figure \ref{f9}. The memory also changes the waiting time.

\section{Final comments}

In this paper we applied the finite element method with polinomial basis of degree $r$ complemented with the Crank-Nicolson method and the trapezoid quadrature for a class of evolution differential equations with $p$-Laplacian and memory. We present a simple and robust numerical method that appears to have an optimal convergence order with which asymptotic behavior and location properties can be observed. This study is a complement to the theoretical part of Antonsev et al. and Almeida et al.

\section*{Acknowledgements}

This work was partially supported by the research projects: Grant N.\linebreak UID/MAT/00212/2019 - financed by FEDER through the - Programa Operacional
Factores de Competitividade, FCT - Funda\c{c}\~{a}o para a Ci\^{e}ncia e a Tecnologia and Grant BID/ICI-FC/Santander Universidades-UBI/2015.


\end{document}